\documentclass[journal,twoside,web]{ieeecolor}
\usepackage{generic}
\usepackage{cite}
\usepackage{amsmath,amssymb,amsfonts}
\usepackage{algorithmic}
\usepackage{graphicx}
\usepackage{algorithm,algorithmic}
\usepackage{hyperref}
\hypersetup{hidelinks=true}
\usepackage{textcomp}
\usepackage{tikz}
\usepackage{epstopdf}

\newtheorem{theorem}{Theorem}[section]

\newtheorem{lemma}[theorem]{Lemma}
\newtheorem{proposition}[theorem]{Proposition}

\newtheorem{remark}[theorem]{Remark}
\newtheorem{assumption}[theorem]{Assumption}

\def\BibTeX{{\rm B\kern-.05em{\sc i\kern-.025em b}\kern-.08em
    T\kern-.1667em\lower.7ex\hbox{E}\kern-.125emX}}
\begin{document}
\title{Stabilization of a Class of Large-Scale Systems of Linear Hyperbolic PDEs via Continuum 
	Approximation of Exact Backstepping Kernels$^*$}
\author{Jukka-Pekka Humaloja, \IEEEmembership{Member, IEEE}, and Nikolaos Bekiaris-Liberis, 
\IEEEmembership{Senior Member, IEEE}
\thanks{$^*$Funded by the European Union (ERC, C-NORA, 101088147). Views and opinions 
	expressed are however those of the authors only and do not necessarily reflect those of the 
	European Union or the European Research Council Executive Agency. Neither the European 
	Union nor the granting authority can be held responsible for them.}
\thanks{$^1$The authors are with the Department of Electrical and Computer Engineering, 
	Technical University of Crete, Chania, Greece. Emails: jhumaloja@tuc.gr and nlimperis@tuc.gr.}
	}

\maketitle

\begin{abstract}
We establish that stabilization of a class of linear, hyperbolic partial differential equations (PDEs) 
with a large (nevertheless 
finite) number of components, can be achieved via employment of a backstepping-based control 
law, which is constructed for stabilization of a continuum version (i.e., as the number of 
components tends to infinity) of the PDE system. This is achieved by proving that the exact 
backstepping kernels, constructed for stabilization of the large-scale system, can be 
approximated 
(in certain sense such that exponential stability is preserved) by the backstepping kernels 
constructed for stabilization of a continuum version (essentially an infinite ensemble) of the 
original PDE system. The proof relies on construction of a convergent sequence of backstepping 
kernels that is defined such that each kernel matches the exact backstepping kernels (derived 
based on the original, large-scale system), in a piecewise constant manner with respect to an 
ensemble variable; while showing that they satisfy the continuum backstepping kernel equations. 
We present a numerical example that reveals that complexity of computation of stabilizing 
backstepping kernels may not scale with the number of components of the PDE state, when the 
kernels are constructed on the basis of the continuum version, in contrast to the case in which 
they are constructed on the basis of the original, large-scale system. In addition, we formally 
establish the connection between the solutions to the large-scale system and its continuum 
counterpart. Thus, this approach can be 
useful for design of computationally tractable, stabilizing backstepping-based control laws for 
large-scale PDE systems. 
\end{abstract}


\section{Introduction}
\label{sec:introduction}

\subsection{Motivation}

\IEEEPARstart{L}{arge-scale}  systems of 1-D hyperbolic PDEs appear in a variety of applications 
involving 
transport phenomena and which incorporate different, interconnected components. Among them, 
large-scale interconnected hyperbolic systems may be used to describe the dynamics of blood 
flow, from the location of the heart all the way through to points where non-invasive 
measurements can be obtained \cite{BikPhd}, \cite{ReyMer09}, of epidemics spreading, to 
describe transport of 
epidemics among different geographical regions \cite{GuaPri2}, and of traffic flows, to model 
density and 
speed dynamics in interconnected highway segments \cite{FriGot22,ZhaLua22} and in urban 
networks 
\cite{TumCan22}, to name a 
few \cite{BasCorBook}. Backstepping is a systematic design approach to construction of explicit 
feedback laws for general 
classes of such systems \cite{AurBre22, KrsSmyBook, Aur20, 
	HuLDiM16,RedAur21,DiMVaz13, EndGab24}. Due to potentially 
large number of 
interacting components, incorporated in such systems, computational complexity of exact 
backstepping-based control designs may increase significantly, in a manner proportional with the 
number of state components. Motivated by this, in the present paper we aim at developing an 
approach to computing backstepping kernels for large-scale hyperbolic PDE systems, such that 
computational complexity remains tractable, even when the number of state components 
becomes very 
large, while at the same time, provably retaining the stability guarantees of backstepping. We 
achieve this via approximating the exact backstepping kernels, computed based on a large 
number of 
PDEs, utilizing a single kernel that is derived based on a continuum version of the exact kernels 
PDEs, and capitalizing on the robustness properties\footnote{Such robustness properties have 
	been also reported within the framework of, e.g., robust output regulation for abstract, 
	infinite-dimensional systems  \cite{PauHumCDC22}. Thus, the approach presented here could 
	be, in principle, also combined with other types of stabilizing feedback laws.} of backstepping 
	for the class of systems considered, to additive control gain errors.

\subsection{Literature}

The approach of design of feedback laws for large-scale systems based on a continuum version 
of the system considered has been utilized for large-scale ordinary differential equation (ODE) 
systems, such as, for example, in 
\cite{BloHen09, FerBuf06, MeuKrs11, NikCan22, WeiFri19, QiJVaz15, FriKrs11, ZhaVaz24}. 
However, such an approach has not been utilized so far for 
large-scale systems whose state components are PDEs. The main goals of the approach 
developed here, may be viewed also as related to control design approaches that aim at providing 
computational means towards implementation of PDE backstepping-based control laws with 
provable stability guarantees, such as, for example, neural operators-based \cite{BhaShi24}, 
late-lumping-based \cite{AurMor19}, and power series-based \cite{VazCheCDC23}, backstepping 
control laws. Our approach 
may be viewed as complementary and different from these results, in that the main goal is to 
address complexity due to a potential radical increase in the number of state components, instead 
of complexity of actual numerical implementation (even though these existing results can be 
combined with the approach presented here, for numerical implementation of the controllers). 

\subsection{Contributions}

In the present paper, we provide backstepping-based feedback laws for a class of large-scale 
systems of 1-D hyperbolic PDEs, which are described by the class of systems considered in 
\cite{DiMVaz13}, 
when the number of state components $n$ is large. The key idea of our approach is to construct 
approximate 
backstepping kernels for stabilization of the large-scale (nevertheless, with a finite number of 
components) system relying on the continuum backstepping kernels developed in 
\cite{AllKrs23arxiv} for a 
continuum version of the original, large-scale system. We establish stability of the closed-loop 
system consisting of the original, large-scale PDE system under a backstepping-based feedback 
law that employs the approximate kernels, constructed based on the continuum version of the 
PDE system. The stability proof consists of three main steps. 

In the first, we construct a sequence of 
backstepping kernels that is defined such that each kernel matches with the exact backstepping 
kernel (derived based on the original, large-scale system) in a piecewise constant manner with 
respect to an ensemble variable; while showing that the kernels in the sequence satisfy the 
continuum kernel equation. 
For the proof we rely on a transformation that maps functions (the exact backstepping kernels) on 
a 2-D domain, to functions (the approximate kernels derived based on the continuum) on a 3-D 
domain in a piecewise constant manner in $L^2$ sense. In the second step, we show that this 
sequence converges to the continuum 
backstepping kernel, obtained from a direct application of backstepping to a continuum version of 
the large-scale system. For the proof we rely on the well-posedness of the backstepping kernels 
(both the continuum and exact kernels) and density arguments. In the third step, we establish 
stability of the closed-loop system employing an abstract systems framework. For the proof we 
recast the closed-loop system’s dynamics as perturbed dynamics of the nominal (based 
on employment of the exact kernels) closed-loop system, showing that the size of the 
perturbation, due to the additive error that originates from the approximation error of the exact 
stabilizing control gains, can be made arbitrarily small (in $L^2$) for sufficiently large $n$. 

We also provide an alternative stability proof employing a Lyapunov functional, which allows 
quantification of overshoot and decay rate of the closed-loop system's response. Furthermore, 
for enabling generalization of the approach introduced, for computation of stabilizing control 
gains based on continuum approximations, to other classes of large-scale PDE systems, we also 
establish the formal connection between the solutions to the original $n+1$ system and the 
solutions to its continuum counterpart. In particular, we show that when the number of state 
components is sufficiently large, the solutions to the large-scale system can be approximated by 
the solutions to the continuum system, provided that the data (i.e., parameters, initial conditions, 
and inputs) of the $n+1$ PDE problem can be approximated by the respective data of the 
continuum PDE problem. The proof relies on construction of a sequence of solutions, obtained in 
a piecewise constant manner (with respect to an ensemble variable) from the solutions to the 
$n+1$ system, which is subsequently shown (via utilization of the well-posedness property of the 
continuum that we prove) to converge to the solution of the continuum system.

We then present a numerical example that illustrates that computation of (approximate) stabilizing 
kernels based on the continuum kernel may provide flexibility in computation, as well as it may 
significantly improve computational complexity. In particular, in this specific example, although 
computation of the exact backstepping kernels may require to solve implicitly the corresponding 
$n+1$ 
hyperbolic kernels PDEs, as closed-form solutions may not be available, the approximate 
kernels can be computed with only algebraic computations, since the continuum kernel is 
available in closed form. We also present respective simulation investigations, which validate the 
theoretical developments, showing that as the number of components of the large-scale system 
increases the performance of the closed-loop systems, under the approximate control laws, is 
improved.  In particular, we illustrate that the approximate control kernels converge to the exact 
kernels, and thus, as $n$ increases, the performance of the closed-loop system becomes similar 
to the performance under the exact control gain kernels.

\subsection{Organization}

We start in Sections~\ref{sec:n+1} and \ref{sec:inf} presenting both the large-scale PDE system 
and its continuum 
counterpart, together with the respective exact and continuum backstepping kernels PDEs. In 
Section~\ref{sec:stab} we establish stability of the large-scale, closed-loop system under the 
approximate control law.  In Section~\ref{sec:sim} we present a 
numerical example and consistent simulation results. In Section~\ref{sec:part} we study the 
connection between the solutions to the large-scale system and its continuum counterpart. In 
Section~\ref{sec:ntoinf} we discuss separately the case $n\to\infty$. In 
Section~\ref{sec:conc} we provide concluding remarks and discuss related topics of our current 
research.

\subsection{Notation}

We use the standard notation $L^2(\Omega; \mathbb{R})$ for real-valued
Lebesque integrable functions on a domain $\Omega$, and on one
dimensional domains $H^1$ denotes the corresponding Sobolev space. Similarly,
$L^\infty(\Omega;\mathbb{R}),C(\Omega; \mathbb{R}), C^1(\Omega; \mathbb{R})$ denote 
essentially bounded, continuous, and
continuously differentiable functions, respectively, on $\Omega$. Moreover,
$f \in L_{\rm{loc}}^2([0,+\infty);\mathbb{R})$ means that $f \in L^2([0,
m];\mathbb{R})$ for any $m\in \mathbb{N}$. We 
denote vectors and matrices by bold symbols, and $\|\cdot\|_\infty$ denotes the maximum 
absolute row sum of a matrix (or a vector). For any $n \in \mathbb{N}$, we denote by $E$ the 
Hilbert space 
$L^2([0,1]; \mathbb{R}^{n+1})$ equipped with the inner product 
\begin{equation}
	\label{eq:eip}
	\left\langle \left( 
	\begin{smallmatrix}
		\mathbf{u}_1 \\ v_1
	\end{smallmatrix}
	\right), \left( 
	\begin{smallmatrix}
		\mathbf{u}_2\\ v_2
	\end{smallmatrix}
	\right) \right\rangle_E
	= \int\limits_0^1 \left(\frac{1}{n}
	\sum_{i=1}^nu_1^i(x)u_2^i(x) + v_1(x)v_2(x)\right)dx,
\end{equation}
which induces the norm $\left\|
\cdot
\right\|_E = \sqrt{\left<\cdot,\cdot\right>_E}$. We also define the continuum version of $E$ as 
$E_c = L^2([0,1];L^2([0,1];\mathbb{R}))\times L^2([0,1];\mathbb{R})$, (i.e., $\mathbb{R}^n$ 
becomes $L^2([0,1];\mathbb{R})$ as $n\to\infty$) equipped with the inner 
product
\begin{align}
	\label{eq:ecip}
	\left\langle \left( 
	\begin{smallmatrix}
		u_1 \\ v_1
	\end{smallmatrix}
	\right), \left( 
	\begin{smallmatrix}
		u_2\\ v_2
	\end{smallmatrix}
	\right) \right\rangle_{E_c}
	& = \nonumber \\
	 \int\limits_0^1\left(\int\limits_0^1 u_1(x,y)u_2(x,y)dy + v_1(x)v_2(x)\right)dx, &
\end{align}
which coincides with $L^2([0,1]^2;\mathbb{R})\times L^2([0,1];\mathbb{R})$.
Moreover, $\mathcal{L}(E, \mathbb{R})$ denotes 
the  space of bounded linear operators from $E$ to $\mathbb{R}$, and 
$\|\cdot\|_{\mathcal{L}(E,\mathbb{R})}$ is the corresponding operator norm. For 
$\mathcal{L}(E,E)$, we denote $\mathcal{L}(E)$. Finally, we say that a 
system is exponentially stable (on $E$; resp. on $E_c$) if for any initial condition $z_0\in E$ (resp. 
$z_0\in E_c$) the (weak) solution 
$z(t)$ of the system satisfies $\|z(t)\|_E \leq Me^{-ct}\|z_0\|_E$ (resp. $\|z(t)\|_{E_c} \leq 
Me^{-ct}\|z_0\|_{E_c}$) for some $M,c > 0$.

\section{Stabilization of Large-Scale Systems of Linear Hyperbolic PDEs via Exact 
	Backstepping Kernels} \label{sec:n+1}

For $n \geq 1$ consider the following set of $n+1$ transport PDEs on $x \in [0,1]$ for
$i = 1,\ldots, n$
\begin{subequations}
	\label{eq:n+1}%
	\begin{align}
		u_t^i(t,x) + \lambda_i(x)u_x^i(t,x)  & = \frac{1}{n}\sum_{j=1}^n\sigma_{i,j}(x)u^j(t,x)  
		\nonumber \\
		& \quad +  W_i(x)v(t,x), \\
		v_t(t,x) - \mu(x)v_x(t,x)
		& = \frac{1}{n}\sum_{j=1}^n \theta_j(x)u^j(t,x), 
	\end{align}
\end{subequations}
with boundary conditions
\begin{equation}
	\label{eq:nbcuy}
	u^i(t,0)  = q_iv(t,0), \qquad v(t,1) = U(t),
\end{equation}
where $U\in L_{\rm{loc}}^2([0, +\infty);\mathbb{R})$ is the control input. The initial conditions of 
\eqref{eq:n+1} are
$u^i(0,x)=u^i_0(x), v(0,x)=v_0(x)$, where $u_0^i, v_0 \in L^2([0,1];\mathbb{R})$. The 
parameters of the system \eqref{eq:n+1}, \eqref{eq:nbcuy} satisfy the following assumption.
\begin{assumption}
	\label{ass:n+1}
	We assume that $\mu, \lambda_i \in C^1([0,1];\mathbb{R})$, 
	$\sigma_{i,j}, W_i,\theta_i \in C([0,1];\mathbb{R})$ and $q_i\in\mathbb{R}$ for all 
	$i,j=1,2,\dots,n$. Moreover, the transport 
	velocities are assumed to satisfy $-\mu(x) < 0 < \lambda_i(x)$,	for all $x \in [0,1]$ and 
	$i=1,2,\ldots,n$.
\end{assumption}
\begin{remark}
	\label{rem:n+1}
	The presentation of the system \eqref{eq:n+1}, \eqref{eq:nbcuy} is motivated from 
	\cite{DiMVaz13}. However, here we also make the following modifications. Most notably, the 
	factor $1/n$ appears in \eqref{eq:n+1}. This is equivalent to equipping
	the $n$-part of the system with the scaled inner product $n^{-1} \langle \cdot, 
	\cdot \rangle_{\mathbb{R}^n}$. With the scaling, we guarantee that the sums remain 
	bounded and convergent as $n\to\infty$ without having to pose any additional constraints on 
	the parameters of \eqref{eq:n+1}, \eqref{eq:nbcuy}. If one wishes to proceed without scaling the 
	sums, then some additional assumptions are needed, e.g., that the respective parameters form 
	$\ell^p$ sequences for some $p\in[1,+\infty]$, such that the sums are well-defined as 
	$n\to\infty$. The other 
	noteworthy modification has to do with Assumption~\ref{ass:n+1}. In \cite{DiMVaz13},
	the transport velocities are assumed to satisfy $-\mu(x) < 0 
	< \lambda_1(x)  < \lambda_2(x) < \ldots < \lambda_n(x), \forall x \in[0,1]$ to guarantee strict 
	hyperbolicity and well-posedness of 
	\eqref{eq:n+1}, \eqref{eq:nbcuy}\footnote{ In specific cases, such an assumption may not be 
		required, for example, for constructing a Lyapunov functional, see, e.g., \cite{AtaDas23, 
			TerAnd20}.}. Nevertheless, well-posedness can be guaranteed under 
	Assumption~\ref{ass:n+1}, e.g., based on \cite[Sect. 13.2]{JacZwaBook} as we show in 
	Proposition~\ref{prop:wp} in Appendix \ref{app:tech}. Due to
	well-posedness, the system \eqref{eq:n+1}, \eqref{eq:nbcuy} has a well-defined, unique, weak 
	solution $\left(\left(u^i(t,x)\right)_{i=1}^n,v(t,x)\right)$ on $E$, where $u^i,v \in 
	C([0,+\infty); L^2([0,1];\mathbb{R}))$ \cite[Prop. 4.2.5, Rem. 4.1.2]{TucWeiBook}.
	\label{rem:scale}
\end{remark}

It follows from \cite[Thm 3.2]{DiMVaz13} that the system  \eqref{eq:n+1}, \eqref{eq:nbcuy} is 
exponentially stabilizable by a state feedback law of the form 
\begin{equation}
	\label{eq:Un}
	U(t) = \int\limits_0^1 \left[  \frac{1}{n}\sum_{i=1}^nk^i(1,\xi)u^i(t,\xi) +
	k^{n+1}(1,\xi)v(t,\xi) \right] d\xi,
\end{equation}
where, for $i=1,\ldots,n+1$, $k^i$ satisfy 
\begin{subequations}
	\label{eq:kn}%
	\begin{align}
		\mu(x)k_x^i(x,\xi) - \lambda_{i}(\xi)k_{\xi}^i(x,\xi) & = \nonumber  \\
		\lambda_i'(\xi)k^i(x,\xi)+\frac{1}{n}\sum_{j=1}^n \sigma_{j,i}(\xi)k^j(x,\xi) +
		\theta_i(\xi)k^{n+1}(x,\xi),& \\
		\mu(x)k_x^{n+1}(x,\xi) + \mu(\xi)k_{\xi}^{n+1}(x,\xi) & = \nonumber \\
		-\mu'(\xi)k^{n+1}(x,\xi) +  \frac{1}{n}\sum_{j=1}^nW_j(\xi)k^j(x,\xi), & 
	\end{align}
\end{subequations}
on a triangular domain $0 \leq \xi \leq x \leq 1$ with boundary conditions
\begin{subequations}
	\label{eq:knbc}%
	\begin{align}
		k^i(x,x) & = -\frac{\theta_i(x)}{\lambda_i(x) + \mu(x)}, \\
		\mu(0)k^{n+1}(x,0) & =  \frac{1}{n}\sum_{j=1}^nq_j\lambda_j(0)k^j(x,0),
	\end{align}
\end{subequations}
for all $x \in [0,1]$. Note the scaling of the sums by $1/n$ as per Remark~\ref{rem:scale}.

\section{Stabilization of a Continuum of Linear Hyperbolic PDEs via Continuum 
	Backstepping} \label{sec:inf}

\begin{figure}[t]
	\begin{center}
		\begin{tikzpicture}
			\filldraw[cyan!10] (1,2) rectangle (7,4);
			\draw[-stealth, thick] (1,4) -- (7,4);
			\node at (4,4.25) {$u(t,x,y)$};
			\draw[-stealth, thick] (1,2) -- (7,2);
			\draw[-stealth, thick] (7,1) -- (1,1);
			\node at (4,1.25) {$v(t,x)$};
			\node[color=magenta] at (7.5,1) {$U(t)$};
			\draw[-, red] (1,2) -- (1,4);
			\draw[-stealth, magenta] (1,1) .. controls (0,2) .. (1,3);
			\node[color=magenta] at (0.65, 2) {$q(y)$};
			\draw[-stealth, blue] (2,1) -- (2,2);
			\node[color=blue] at (2.6,1.5) {$W(x,y)$};
			\draw[-stealth, blue] (5,2) -- (5,1);
			\node[color=blue] at (5.95,1.5) {$\displaystyle \int\theta(x,y)\cdot dy$};
			\node[color=blue] at (4,3) {$\displaystyle \int \sigma(x,y,\eta)\cdot d\eta$};
			\draw[-stealth] (7.5, 2) -- (7.5, 4);
			\node at (7.75,4) {$y$};
			\draw[-stealth] (1,0.5) -- (7, 0.5);
			\node at (7,.25) {$x$};
		\end{tikzpicture}
	\end{center}
	\caption{Schematic view of the continuum PDE system \eqref{eq:inf}, \eqref{eq:cbcuy}. 
	Boundary 
		terms are denoted in magenta and in-domain terms are denoted in blue.}
	\label{fig:inf}
\end{figure}
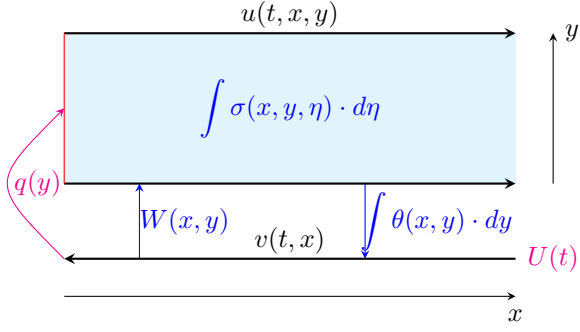

While large-scale, yet, consisting of a finite-number of components, systems of hyperbolic PDEs 
can be studied in the framework of Section 
\ref{sec:n+1}, we also consider the continuum limit case as $n\to\infty$, for which we 
present the 
generic framework of a continuum of hyperbolic PDEs studied in \cite{AllKrs23arxiv} and sketched 
in Fig.~\ref{fig:inf}. That is, 
instead of having $n$ rightward transport PDEs as in Section~\ref{sec:n+1}, consider a continuum 
of such PDEs as in \cite{AllKrs23arxiv} with $y\in[0,1]$ being the index variable\footnote{Note that 
	we have not yet formally proved that the continuum limit of system \eqref{eq:n+1} as 
	$n\to\infty$ 
	is system \eqref{eq:inf}. However, we use \eqref{eq:inf} here as an educated guess of such a 
	continuum version (see also \cite{AllKrs23arxiv}) to obtain a continuum version of the 
	respective backstepping kernels $k^i, i=1,\ldots, n$, given in \eqref{eq:kn}. In 
	Section~\ref{sec:part} we, in fact, formally prove that \eqref{eq:inf}, \eqref{eq:cbcuy} is the 
	continuum limit of \eqref{eq:n+1}, \eqref{eq:nbcuy}.}
\begin{subequations}
	\label{eq:inf}%
	\begin{align}
		u_t(t,x,y) + \lambda(x,y)u_x(t,x,y) \nonumber  & = \int\limits_0^1 
		\sigma(x,y,\eta)u(t,x,\eta)d\eta 
		\nonumber \\
		& \quad  + W(x,y)v(t,x),   \\
		v_t(t,x) - \mu(x)v_x(t,x)  & = \int\limits_0^1\theta(x,y)u(t,x,y)dy,
	\end{align}
\end{subequations}
with boundary conditions 
\begin{equation}
	\label{eq:cbcuy}
	u(t,0,y) = q(y)v(t,0), \qquad v(t,1) = U(t),
\end{equation}
for almost every $y \in [0,1]$, i.e., the continuum variables and parameters are considered 
$L^2([0,1];\mathbb{R})$ functions in $y$. The following assumption is needed, for the parameters 
involved in \eqref{eq:inf}, \eqref{eq:cbcuy}, to guarantee the existence of a continuum 
backstepping control law \cite[Thm 3]{AllKrs23arxiv}.

\begin{assumption}
	We assume that $\mu \in C^1([0, 1]; \mathbb{R}), \lambda\in C^1\left([0,1]^2;\mathbb{R}\right), 
	W,\theta\in C([0,1];L^2([0,1];\mathbb{R})), \sigma \in 
	C\left([0,1];L^2([0,1]^2;\mathbb{R})\right)$, and $q \in C([0,1];\mathbb{R})$. 
	Moreover, $\lambda(x,y) > 0$ uniformly for all $x\in [0,1]$ and almost every $y\in[0,1]$, 
	and $-\mu(x) < 0$ for all $x\in[0,1]$.
	\label{ass:inf}
\end{assumption}

By \cite[Thm 1]{AllKrs23arxiv}, the system \eqref{eq:inf}, \eqref{eq:cbcuy} is exponentially 
stabilizable by a state feedback law of the form 
\begin{equation}
	\label{eq:Uc}
	U(t) = \int\limits_0^1 \left[\int\limits_0^1 k(1,\xi,y)u(t,\xi,y)dy
	+ \bar{k}(1,\xi)v(t,\xi) \right]d\xi,
\end{equation}
where $k, \bar{k}$ satisfy 
\begin{subequations}
	\label{eq:kc}%
	\begin{align}
		\mu(x)k_x(x,\xi,y) - \lambda(\xi,y)k_{\xi}(x,\xi,y) -  \theta(\xi,y)\bar{k}(x,\xi) & = \nonumber \\
		\lambda_{\xi}(\xi,y)k(x,\xi,y) 
		+\int\limits_0^1\sigma(\xi,\eta,y)k(x,\xi,\eta)d\eta, & \\
		\mu(x) \bar{k}_x(x,\xi) + \mu(\xi)\bar{k}_{\xi}(x,\xi) & =  \nonumber \\
		-\mu'(\xi)\bar{k}(x,\xi) + \int\limits_0^1W(\xi,y)k(x,\xi,y)dy, &
	\end{align}
\end{subequations}
on a triangular domain $0 \leq \xi \leq x \leq 1$ with boundary conditions 
\begin{subequations}
	\label{eq:kcbc}%
	\begin{align}
		k(x,x,y) & = - \frac{\theta(x,y)}{\lambda(x,y) + \mu(x)}, \\
		\mu(0)\bar{k}(x,0) & = \int\limits_0^1q(y)\lambda(0,y)k(x,0,y)dy,
	\end{align}
\end{subequations}
for almost every $y\in [0,1]$.

\section{Stabilization of the Finite Large-Scale System via Continuum Approximation of 
	Exact Kernels} \label{sec:stab}

\subsection{Statement of the Main Result}

The core idea of the continuum approximation  that we present here is that we approximate $n+1$ 
kernel equations by a continuum of kernel equations. Provided that the approximation is 
sufficiently accurate, we show that the backstepping controller derived from the continuum 
kernel equations exponentially stabilizes the $n+1$ system associated with the original finite 
system of $n+1$ PDEs.

Thus, consider an $n+1$ system  \eqref{eq:n+1}, \eqref{eq:nbcuy} with parameters $\lambda_i, 
W_i, \theta_i, \sigma_{i,j}, q_i$ for $i,j=1,2,\ldots,n$ satisfying 
Assumption~\ref{ass:n+1} and consider any continuous functions $\lambda, W, \theta, \sigma,q$ 
that satisfy Assumption~\ref{ass:inf} with
\begin{subequations}
	\label{eq:afn}%
	\begin{align}
		\lambda(x,i/n) & = \lambda_i(x), \\
		W(x, i/n) & = W_i(x), \\
		\theta(x,i/n) & = \theta_i(x), \\
		\sigma(x,i/n,j/n) & = \sigma_{i,j}(x), \label{eq:afns}\\	
		q(i/n) & = q_i,
	\end{align}%
\end{subequations}
for all $x\in[0,1]$ and $i,j = 1,2,\ldots,n$. There are infinitely many functions 
satisfying \eqref{eq:afn} and Assumption~\ref{ass:inf}, as well as ways to construct them, e.g., by 
utilizing auxiliary functions $\left\{p_i\right\}_{i=1}^n$ that satisfy $p_i(i/n)=1$ and 
$p_\ell(i/n)=0$ for $\ell\neq i$.\footnote{\label{fn:afn}We demonstrate this by constructing  
	$\sigma(x,y,\eta)$ satisfying  \eqref{eq:afns} of the form
	$\sigma(x,y,\eta) =  \sum_{i=1}^n\sum_{j=1}^n \sigma_{i,j}(x)p_i(y)p_j(\eta)$, where
	\begin{equation}
		p_i(y) = \prod_{k=1,k\neq i}^{n} \frac{(k/n - y)}{(k/n - i/n)} + b\sin(n\pi y),
	\end{equation}
	satisfies $p_i(i/n) = 1$, $p_\ell(i/n) = 0$ for $\ell\neq i$, and 
	$p_i\in C^\infty([0,1];\mathbb{R})$ for any $i = ,1,2,\ldots,n$ and any $b\in\mathbb{R}$. 
	Similar constructions can be obtained for $\lambda, W, \theta$, and $q$.} 
The relations in \eqref{eq:afn} could as well be defined in other ways, e.g., using $(i-1)/n$ in 
place of $i/n$, but we find \eqref{eq:afn} the most convenient option for our developments. The 
continuum kernel 
equations \eqref{eq:kc}, \eqref{eq:kcbc} with parameters $\lambda, \mu,
W, \theta, \sigma,q$, satisfying \eqref{eq:afn} and Assumption~\ref{ass:inf} have a unique, 
continuous solution $(k, \bar{k})$ by \cite[Thm 
3]{AllKrs23arxiv}. Thus, construct the following functions for all $0 \leq \xi \leq x \leq 1$
\begin{subequations}
	\label{eq:akn}%
	\begin{align}
		\widetilde{k}^i(x,\xi) & = k(x,\xi, i/n), \quad i = 1,2,\ldots,n, \\
		\widetilde{k}^{n+1}(x,\xi) & = \bar{k}(x,\xi).
	\end{align}
\end{subequations}
Our main result is the following.
\begin{theorem}
	\label{thm:stab}
	Consider an $n+1$ system \eqref{eq:n+1}, \eqref{eq:nbcuy} with parameters 
	$\lambda_i,\mu,W_i,\theta_i,\sigma_{i,j},q$ for $i,j=1,2,\ldots,n$ satisfying
	Assumption~\ref{ass:n+1}. Let the parameters $\lambda,\mu, W, \theta, \sigma, q$ satisfy
	Assumption~\ref{ass:inf} and relations \eqref{eq:afn}. Then, if $n$ is sufficiently large, the 
	system \eqref{eq:n+1}, \eqref{eq:nbcuy} is exponentially stabilized by the control law
	\begin{align}
		\label{eq:aUn}
		U(t) & = \int\limits_0^1 \left[  \frac{1}{n}\sum_{i=1}^n\widetilde{k}^i(1,\xi)u^i(t,\xi) +
		\widetilde{k}^{n+1}(1,\xi)v(t,\xi) \right]d\xi,
	\end{align}
	where $\left(\widetilde k^i\right)_{i=1}^{n+1}$ are given in \eqref{eq:akn}, with 
	$\left(k,\bar{k}\right)$ being the solution to \eqref{eq:kc}, \eqref{eq:kcbc}.
\end{theorem}

\subsection{Proof of Theorem \ref{thm:stab}}

The proof of Theorem~\ref{thm:stab} relies on Lemmas~\ref{lem:se} and~\ref{lem:ka} presented 
below. We show first that the functions defined in \eqref{eq:akn} approximate the solutions to the 
$n+1$ kernel equations \eqref{eq:kn}, \eqref{eq:knbc} to arbitrary accuracy as $n$ increases. In 
order to do this, we first interpret the solutions to the $n+1$ kernels equations 
\eqref{eq:kn}, \eqref{eq:knbc} as piecewise constant solutions with respect to $y$, to the 
continuum kernels equations \eqref{eq:kc}, \eqref{eq:kcbc}. One way to do this is 
highlighted in the following lemma, which is to transform the $\mathbb{R}^n$-valued 
components of the $n+1$ kernel equations \eqref{eq:kn}, \eqref{eq:knbc} into step 
functions in $y\in[0,1]$.
\begin{lemma}
	\label{lem:se}
	Consider the $n+1$ kernel equations \eqref{eq:kn}, \eqref{eq:knbc} where the parameters 
	satisfy Assumption~\ref{ass:n+1} and  define the following 
	functions for all $x\in[0,1]$, piecewise in $y$ for $i,j = 1,2,\ldots,n$\footnote{These functions 
		can be extended to $y\in[0,1]$ by assigning the value at $y=0$ arbitrarily, which does not 
		affect 
		the functions in the $L^2$ sense. The same applies to \eqref{eq:ensk}.}
	\begin{subequations}
		\label{eq:sfn}%
		\begin{align}
			\lambda^n(x,y) & = \lambda_i(x), \quad y \in
			(i-1)/n, i/n], \\
			\sigma^n(x,y,\eta) & = \sigma_{i,j}(x), \quad y \in
			((i-1)/n,i/n],\\
			& \qquad \qquad \qquad \eta \in ((j-1)/n, j/n], \\
			W^n(x,y) & = W_i(x), \quad y \in ((i-1)/n, i/n], \\
			\theta^n(x,y) & = \theta_i(x), \quad y \in
			((i-1)/n,i/n], \\
			q^n(y) & = q_i, \quad y \in ((i-1)/n, i/n].
		\end{align}
	\end{subequations}
	Construct the following function for all 
	$0\leq \xi \leq x \leq 1$, piecewise in $y$ for $i=1,2,\ldots,n$
	\begin{equation}
		k^n(x,\xi,y) = k^i(x,\xi), \quad y \in ((i-1)/n, i/n],
		\label{eq:ensk}
	\end{equation}
	where $\left(k^i\right)_{i=1}^{n+1}$ is the solution to  \eqref{eq:kn}, \eqref{eq:knbc}. Then, 
	$\left(k^n(x,\xi,y), k^{n+1}(x,\xi)\right)$ satisfies the kernel equations \eqref{eq:kc}, 
	\eqref{eq:kcbc}
	for the parameters defined in \eqref{eq:sfn} and the original $\mu$.

	\begin{proof}
		The claim follows after applying a linear transform to the  
		kernel equations \eqref{eq:kn}, \eqref{eq:knbc}. In order to rigorously present the 
		transformation, we have to write \eqref{eq:kn} as a single equation on 
		$\mathbb{R}^{n+1}$. We introduce the following notation
		\begin{subequations}
			\label{eq:ensnk}%
			\begin{align}
				\mathbf{k}(x,\xi) & = \left(k^i(x,\xi)\right)_{i=1}^{n+1}, \\
				\mathbf{L}(x) & = \operatorname{diag}\left(\left(\lambda^i(x)\right)_{i=1}^n,
				-\mu(x)\right), \\
				\mathbf{S}(x) & = \frac{1}{n}
				\begin{bmatrix}
					\sigma_{1,1}(x) & \cdots & \sigma_{n,1}(x) & n\theta_1(x) \\
					\vdots & \ddots & \vdots & \vdots \\
					\sigma_{1,n}(x) & \cdots & \sigma_{n,n}(x) & n\theta_n(x) \\
					W_1(x) & \cdots & W_n(x) & 0
				\end{bmatrix},                             
			\end{align}
		\end{subequations}
		so that \eqref{eq:kn} can be written as 
		\begin{equation}
			\label{eq:ensko}
			\mu(x)\mathbf{k}_x(x,\xi) - \mathbf{L}(\xi)\mathbf{k}_{\xi}(x,\xi) =
			\mathbf{L}'(\xi)\mathbf{k}(x,\xi) + \mathbf{S}(\xi)\mathbf{k}(x,\xi).
		\end{equation}
		The linear transform is given by $\mathcal{F} =
		\operatorname{diag}(\mathcal{F}_{n}, 1)$, where $\mathcal{F}_n\mathbf{e}_{\ell}
		= \chi_{((\ell-1)/n,\ell/n]}$ with  $\chi_{((\ell-1)/n,\ell/n]}$ being the indicator function of the 
		interval 
		$((\ell-1)/n,\ell/n]$ and $\left(\mathbf{e}_{\ell}\right)_{\ell=1}^n$ being
		the Euclidean basis of $\mathbb{R}^n$. Thus, the transform maps any $\mathbf{b}  
		=\left(b_i\right)_{i=1}^{n+1} \in \mathbb{R}^{n+1}$ into $L^2\left([0,1];\mathbb{R}\right)\times 
		\mathbb{R}$ as
		\begin{equation}
			\mathcal{F}\mathbf{b} = \begin{bmatrix}
				\sum_{i=1}^n b_i \chi_{((i-1)/n,i/n]} \\ b_{n+1}.
			\end{bmatrix}.
		\end{equation}
		For any $g \in L^2([0,1];\mathbb{R})$,
		the adjoint $\mathcal{F}_n^{*}$ satisfies
		\begin{equation}
			\left\langle \mathcal{F}_n\mathbf{e}_\ell, g
			\right\rangle_{L^2([0,1];\mathbb{R})}
			= \int\limits_{(\ell-1)/n}^{\ell/n}g(y)dy =
			\frac{1}{n}\left\langle \mathbf{e}_{\ell}, \mathcal{F}_n^{*}g \right\rangle_{\mathbb{R}^{n}},
		\end{equation}
		that is, $\mathcal{F}_n^{*}$ is given by
		\begin{equation}
			\label{eq:Fns}
			\mathcal{F}_n^{*}g = \left( n\int\limits_{(i-1)/n}^{i/n}g(y)dy \right)_{i=1}^n,
		\end{equation}
		where each component is the mean value of $g$ over the interval $[(i-1)/n,i/n]$.
		Thus,  $\mathcal{F}$ has the adjoint
		$\mathcal{F}^{*} = \operatorname{diag}\left(\mathcal{F}_n^{*},1\right)$, which
		additionally satisfies $\mathcal{F}^{*}\mathcal{F} = I_{n+1}$, i.e., $\mathcal{F}$ (and 
		$\mathcal{F}_n$) are isometries, and thus, norm preserving from their domain to their 
		co-domain.
		
		Let us now transform \eqref{eq:ensko} from $\mathbb{R}^{n+1}$ to
		$L^2([0,1];\mathbb{R})\times \mathbb{R}$ by applying $\mathcal{F}$ to
		\eqref{eq:ensko} from the left
		\begin{align}
			\label{eq:enskt}
			\mu(x) \mathcal{F}\mathbf{k}_x(x,\xi) -
			\mathcal{F}\mathbf{L}(\xi)\mathcal{F}^{*}\mathcal{F}\mathbf{k}_{\xi}(x,\xi)
			& = \nonumber \\
			\mathcal{F}\mathbf{L}'(\xi)\mathcal{F}^{*}\mathcal{F}\mathbf{k}(x,\xi)
			+ \mathcal{F}\mathbf{S}(\xi)\mathcal{F}^{*}\mathcal{F}\mathbf{k}(x,\xi),
		\end{align}
		where we also utilized $\mu(x)$ being scalar-valued and
		$\mathcal{F}^{*}\mathcal{F} = I_{n+1}$. Let $g^n$ be a step function in $y$ and $r$ a scalar. 
		Applying the transformations gives 
		\begin{subequations}
			\label{eq:ensft}%
			\begin{align}
				\mathcal{F}\mathbf{k}(x,\xi)
				& = 
				\begin{bmatrix}
					k^n(x,\xi,\cdot) \\ k^{n+1}(x,\xi)
				\end{bmatrix}, \\
				\mathcal{F}\mathbf{L}(\xi)\mathcal{F}^{*}\begin{bmatrix}
					g^n\\ r
				\end{bmatrix}
				& = \begin{bmatrix}
					\lambda^n(\xi,\cdot)g^n(\cdot) \\ \mu(\xi)r
				\end{bmatrix}, \\
				\mathcal{F}\mathbf{S}(\xi)\mathcal{F}^{*}\begin{bmatrix}
					g^n \\ r
				\end{bmatrix}
				& =
				\begin{bmatrix}
					\int\limits_0^1\sigma^n(\xi,\eta,\cdot)g^n(\eta) d\eta & \theta^n(\xi,\cdot)r \\
					\int\limits_0^1W^n(\xi,y)g^n(y) dy  & 0
				\end{bmatrix},
			\end{align}
		\end{subequations}
		with the functions $k^n, \lambda^n, \sigma^n, \theta^n, W^n$ defined
		in \eqref{eq:sfn} and \eqref{eq:ensk}. Inserting \eqref{eq:ensft} into \eqref{eq:enskt} and 
		writing 
		the equations separately on $L^2([0,1];\mathbb{R})$ and $\mathbb{R}$ yields  \eqref{eq:kc} 
		with 
		parameters  $\lambda^n, \mu, \sigma^n, \theta^n, W^n$. Essentially, this amounts to 
		$\left(k^n,k^{n+1}\right)$ satisfying \eqref{eq:kc} on intervals for $y$ of the 
		form \eqref{eq:ensk}, which in turn implies satisfying \eqref{eq:kc} for almost all $y\in[0,1]$, 
		i.e., in the $L^2$ sense with respect to $y$.
		
		The boundary conditions \eqref{eq:knbc} could be transformed into
		\eqref{eq:kcbc} with the same transformation, but it is more
		straightforward to check directly that $\left(k^n,k^{n+1}\right)$
		satisfies \eqref{eq:kcbc} for parameters $\theta^n, \lambda^n, \mu,
		q^n$. For any $i=1,2,\ldots,n$, we have for all $x\in[0,1]$ and $y \in ((i-1)/n,i/n]$
		\begin{equation}
			k^n(x,x,y) = k^i(x,x)
			= \frac{\theta_i(x)}{\lambda_i(x) + \mu(x)}
			=  \frac{\theta^n(x,y)}{\lambda^n(x,y) + \mu(x)},
		\end{equation}
		and thus, this boundary condition is satisfied for almost all
		$y\in[0,1]$. Moreover,
		\begin{align}
			\mu(0)k^{n+1}(x,0)
			& = \frac{1}{n}\sum_{j=1}^nq_j\lambda_j(0)k^j(x,0) \nonumber  \\
			& = \int\limits_0^1q^n(y)\lambda^n(0,y)k^n(x,0,y)dy, 
		\end{align}
		which concludes that the boundary conditions \eqref{eq:kcbc} are
		satisfied. This concludes the proof.
	\end{proof}
\end{lemma}

Let us next consider the continuum of  kernel equations \eqref{eq:kc}, \eqref{eq:kcbc} with 
continuous parameters $\lambda, \mu, W, \theta, \sigma,q$ that satisfy \eqref{eq:afn} 
and Assumption~\ref{ass:inf}, together with the respective kernel equations  \eqref{eq:kc}, 
\eqref{eq:kcbc} with 
piecewise constant parameters $\lambda^n, \mu, W^n, \theta^n, \sigma^n,q^n$ in $y$
constructed in Lemma~\ref{lem:se}. In the next lemma, we show that the solution 
$\left(k^n,\bar{k}^n\right)$ to the latter approximates the solution 
$(k,\bar{k})$ to the former to arbitrary accuracy, provided that $n$ is 
sufficiently 
large.

\begin{lemma}
	\label{lem:ka}
	Consider the solutions $\left(k^n,\bar{k}^n\right)$ to the kernel equations 
	\eqref{eq:kc}, 
	\eqref{eq:kcbc} with parameters $\lambda^n, \mu, W^n,\theta^n, \sigma^n,q^n$ from Lemma 
	\ref{lem:se}. There exist continuous parameters $\lambda, \mu, W, 
	\theta, \sigma,q$ constructed such that they satisfy Assumption~\ref{ass:inf} and 
	\eqref{eq:afn}, 
	and for any such parameters the solution 
	$\left(k,\bar{k}\right)$ to the respective  kernel equations \eqref{eq:kc}, 
	\eqref{eq:kcbc} 
	exists and satisfies the following implications. For any 
	$\delta > 0$, there exists an $n_\delta \in \mathbb{N}$ such that for all $n \geq n_\delta$ we 
	have
	\begin{subequations}
		\label{eq:ka}
		\begin{align}
			\max_{(x,\xi)\in \mathcal{T}}\|k(x,\xi,\cdot) - k^n(x,\xi,\cdot)\|_{L^2([0,1];\mathbb{R})} & \leq 
			\delta, \\
			\max_{(x,\xi)\in \mathcal{T}} |\bar{k}(x,\xi) - \bar{k}^n(x,\xi)| & \leq \delta,
		\end{align}
	\end{subequations}
	where we denote $\mathcal{T} = \{(x,\xi) \in [0,1]^2: 0 \leq \xi \leq x \leq 1\}$.

	\begin{proof}
		We begin by establishing some key properties of the solutions $\left(k^n, \bar{k}^n\right)$ 
		and $(k,\bar k)$.	Firstly, it has been shown in \cite[Sect. V]{DiMVaz13} that, under 
		Assumption~\ref{ass:n+1}, the 
		kernel equations \eqref{eq:kn}, \eqref{eq:knbc} are well-posed, i.e., that the solution 
		$\left(k^i\right)_{i=1}^{n+1}$ exists, is unique, and depends continuously on the 
		parameters of  
		\eqref{eq:kn}, \eqref{eq:knbc}, and that $\left(k^i\right)_{i=1}^{n+1}$ is continuous on 
		$\mathcal{T}$.	Secondly, the functions $\left(k^n, \bar{k}^n\right)$ are 
		constructed in Lemma 
		\ref{lem:ka} 
		based on $\left(k^i\right)_{n=1}^{n+1}$, and thus, for almost all $y$, they exist, are 
		unique, continuous on $\mathcal{T}$, 
		and depend 
		continuously on $\lambda^n, \mu, W^n,\theta^n, \sigma^n,q^n$ in the $L^2$ sense (in $y$). 
		Thirdly, 
		the existence and uniqueness of the solution $\left(k,\bar k\right)$ to the kernel 
		equations \eqref{eq:kc}, \eqref{eq:kcbc} follows, provided that the parameters 
		$\lambda,\mu,W,\theta,\sigma,q$ satisfy Assumption~\ref{ass:inf} \cite[Thm 
		3]{AllKrs23arxiv}. Moreover, it has been shown in \cite[Sect. VI]{AllKrs23arxiv} that 
		$(k,\bar k)$ are continuous on $\mathcal{T}$, and as a consequence of the 
		estimates in \cite[Sect. VI.C]{AllKrs23arxiv}, the solution $\left(k,\bar k\right)$ 
		depends continuously on $\lambda, \mu, W, \theta, \sigma,q$. 
		
		Due to the continuity of $\left(k^n, \bar{k}^n\right)$ and $\left(k,\bar 
		k\right)$  on 
		$\mathcal{T}$, the norms in \eqref{eq:ka} are continuous on $\mathcal{T}$, and hence, the 
		maxima are reached at some point in $\mathcal{T}$. Thus, it remains to show that the 
		maxima 
		become arbitrarily small when $n$ is sufficiently large. The remainder of the proof utilizes the 
		well-posedness of the kernel equations \eqref{eq:kc}, 
		\eqref{eq:kcbc}, in particular that the 
		solutions $\left(k,\bar{k}\right)$ and $\left(k^n,\bar{k}^n\right)$ depend 
		continuously on 
		the parameters of the respective kernel equations as we established in the 
		beginning of the proof. First we show that, when $n$ is sufficiently large, the parameters  
		$\lambda^n, W^n, \theta^n, \sigma^n,q^n$ can approximate the respective continuum 
		parameters  $\lambda, W, \theta, \sigma,q$ to arbitrary accuracy in the $L^2$ sense. That is, 
		for any $\varepsilon_1 > 0$, the following estimates are satisfied for any sufficiently 
		large~$n$
		\begin{subequations}
			\label{eq:sa}%
			\begin{align}
				\max_{x\in[0,1]}\|\lambda(x,\cdot) - \lambda^n(x,\cdot)\|_{L^2([0,1];\mathbb{R})} & \leq 
				\varepsilon_1, \\
				\max_{x\in[0,1]}\|\sigma(x,\cdot,\cdot) - 
				\sigma^n(x,\cdot,\cdot)\|_{L^2([0,1]^2;\mathbb{R})}
				& 
				\leq 	\varepsilon_1, \\
				\max_{x\in[0,1]}\|\theta(x,\cdot) - \theta^n(x,\cdot)\|_{L^2([0,1];\mathbb{R})} & \leq 
				\varepsilon_1, \\
				\max_{x\in[0,1]}\|W(x,\cdot) - W^n(x,\cdot)\|_{L^2([0,1];\mathbb{R})} & \leq 
				\varepsilon_1, \\
				\|q-q^n\|_{L^2([0,1];\mathbb{R})} & \leq 
				\varepsilon_1.
			\end{align}%
		\end{subequations}
		Existence of continuous functions $\lambda, W, \theta, \sigma,q$ satisfying \eqref{eq:afn} 
		and 
		Assumption~\ref{ass:inf} can be guaranteed by construction (see, e.g., 
		Footnote~\ref{fn:afn}). Moreover, the 
		functions  $\lambda^n, W^n, \theta^n, \sigma^n,q^n$ are continuous in $x$ by
		construction and Assumption~\ref{ass:n+1}, so that the maxima in \eqref{eq:sa} are 
		reached at some points $x\in[0,1]$. As the functions compared in \eqref{eq:sa} match 
		for all $x\in[0,1]$ at the points $y,\eta=1/n,2/n,\ldots,1$, the differences of the functions 
		can be made  arbitrarily small in the sense of \eqref{eq:sa} with $n$ being sufficiently large, 
		e.g., based on the fact that step functions are dense in the $L^2$ function space  \cite[Sect. 
		1.3.5]{TaoBook11}.
		
		Finally, due to the well-posedness of the kernel equations, the solutions compared in 
		\eqref{eq:ka} depend continuously on the parameters compared in \eqref{eq:sa}. Thus,  as 
		$\varepsilon_1$ tends to zero in \eqref{eq:sa}, the difference of the solutions 
		$\left(k,\bar{k}\right)$ and 
		$\left(k^n,  \bar{k}^n\right)$ converges to zero in a certain sense. More precisely, the 
		convergence is 
		exactly in the sense stated in \eqref{eq:ka}, as the parameter functions converge uniformly 
		in $x$ by \eqref{eq:sa}, and $\xi$ only appears in place 
		of $x$ in the kernel equations. Thus, the convergence of the solutions is uniform in both $x$ 
		and $\xi$ on $\mathcal{T}$, and hence, for any $\delta > 0$, the estimates \eqref{eq:ka} are 
		satisfied for any sufficiently large~$n$.
	\end{proof}
\end{lemma}

\begin{remark}
	The convergence of the solutions to the kernel equations \eqref{eq:kc}, \eqref{eq:kcbc}, with 
	parameters $\lambda, \mu, W, \theta, \sigma,q$ holds true for any step functions that 
	approximate the continuous parameters as in \eqref{eq:sa} and not only the ones defined in 
	\eqref{eq:sfn}. Any such construction should be such that $\lambda^n, W, \theta^n, 
	\sigma^n,q^n$ match with $\lambda, W, \theta, \sigma,q$ at some points in $y$ within 
	intervals of the form $\left[(i-1)/n, i/n\right]$.
\end{remark}

We have by now established convergence of the solutions $\left(k^n, \bar{k}^n\right)$ to the 
kernel equations \eqref{eq:kc}, \eqref{eq:kcbc} with parameters $\lambda^n,\mu,W^n, \theta^n, 
\sigma^n, q^n$ defined in \eqref{eq:sfn}, to the solutions $\left(k, \bar k\right)$ to \eqref{eq:kc}, 
\eqref{eq:kcbc} with parameters $\lambda, \mu, W, \theta, \sigma,q$. Since the solutions 
$\left(k^n, \bar k^n\right)$ are piecewise constant in $y$ satisfying \eqref{eq:ensk} and 
\eqref{eq:ka}, where 
$\left(k^i\right)_{k=1}^{n+1}$ are 
the solutions to \eqref{eq:kn}, \eqref{eq:knbc} with parameters 
$\left(\lambda^i\right)_{i=1}^n,\mu,\left(W^i\right)_{i=1}^n,\left(\theta\right)_{i=1}^n,\left(q^i\right)_{i=1}^n$,
the solutions $\left(k^n , 
\bar k^n\right)$ can, in fact, approximate the kernels $\left(k^i\right)_{i=1}^{n+1}$ to arbitrary 
accuracy as $n$ 
gets sufficiently large. This in turn implies that the control law \eqref{eq:aUn}, constructed based 
on the solutions $\left(k, \bar k\right)$ to the continuum kernel equations \eqref{eq:kc}, 
\eqref{eq:kcbc}, 
approximates (arbitrarily close as $n$ gets sufficiently large) the original control law \eqref{eq:Un} 
constructed based on the solutions 
$\left(k^i\right)_{i=1}^{n+1}$ to the kernels equations \eqref{eq:kn}, \eqref{eq:knbc}. For this, we 
present the following lemma.

\begin{lemma}
	\label{lem:aU}
	The control law \eqref{eq:aUn} can be written as
	\begin{equation}
		\label{eq:aU}
		\resizebox{.99\columnwidth}{!}{$
			\begin{split}	
				U(t) & = \int\limits_0^1 \left( \frac{1}{n}\sum_{i=1}^nk^i(1,\xi)u^i(t,\xi) +
				k^{n+1}(1,\xi)v(t,\xi) \right) d\xi \\
				& \quad + \int\limits_0^1 \left(  \frac{1}{n}\sum_{i=1}^n\Delta k^i(1,\xi)u^i(t,\xi) +
				\Delta k^{n+1}(1,\xi)v(t,\xi) \right) d\xi,
			\end{split}$}
	\end{equation}
	where $\left(k^i\right)_{i=1}^{n+1}$ is the solution to the $n+1$ kernel equations 
	\eqref{eq:kn}, \eqref{eq:knbc}, and the approximation error terms
	$\left(\Delta k^i(1,\cdot)\right)_{i=1}^{n+1}$ become arbitrarily small, uniformly in $\xi\in[0,1]$, 
	when 
	$n$ 
	is sufficiently large.
	
	\begin{proof}
		Transform the functions $\left(\widetilde{k}^i\right)_{i=1}^{n}$ from \eqref{eq:aUn} into 
		a step function in $y$ as 
		\begin{equation}
			\label{eq:kts}
			\widetilde{k}^n	(x,\xi,y) = \widetilde{k}^i(x,\xi), \quad y \in ((i-1)/n,i/n],
		\end{equation}
		for all $0 \leq \xi \leq x \leq 1$ and piecewise in $y$ for $i=1,2,\ldots,n$. By \eqref{eq:kts} and 
		\eqref{eq:akn}, the functions $\widetilde{k}$ and $k$ coincide for all $0\leq \xi \leq x \leq 1$ 
		at 
		every $y=1/n,2/n,\ldots,1$. As both functions are continuous on $0\leq \xi\leq x \leq 1$ and 
		$k$ is 
		additionally continuous in $y$, for any
		$\varepsilon > 0$, there exists some $n_\varepsilon \in \mathbb{N}$ such that 
		\begin{equation}
			\label{eq:kta}
			\max_{(x,\xi)\in \mathcal{T}}\| k(x,\xi,\cdot) - 	\widetilde{k}^n	
			(x,\xi,\cdot)\|_{L^2([0,1];\mathbb{R})} \leq \varepsilon,
		\end{equation}
		for any $n\geq n_\varepsilon$. Combining this with the estimates of Lemma~\ref{lem:ka} and 
		using 
		the triangle inequality, we have for any $n \geq \max\{n_\delta,n_\varepsilon\}$
		\begin{align}
			\label{eq:ktsa}
			& \quad \max_{(x,\xi)\in \mathcal{T}}\| k^n(x,\xi,\cdot) - 	\widetilde{k}^n	
			(x,\xi,\cdot)\|_{L^2([0,1];\mathbb{R})}  \nonumber \\
			& \leq \max_{(x,\xi)\in \mathcal{T}}\| k^n(x,\xi,\cdot) - 	
			k(x,\xi,\cdot)\|_{L^2([0,1];\mathbb{R})} \nonumber \\
			& \qquad + \max_{(x,\xi)\in \mathcal{T}}\| k(x,\xi,\cdot) - 	
			\widetilde{k}^n	(x,\xi,\cdot)\|_{L^2([0,1];\mathbb{R})}  \nonumber \\
			& \leq \delta + \varepsilon,
		\end{align}
		where both $\delta$ and $\varepsilon$ can be made arbitrarily small by increasing $n$ by 
		Lemma 
		\ref{lem:ka} and \eqref{eq:kta}, respectively. As the estimate is uniform on $\mathcal{T}$, it 
		particularly applies on $x=1$. 
		
		Moreover, the step functions $\widetilde{k}^n$ and $k^n$ constructed in \eqref{eq:kts} 
		and \eqref{eq:ensk}, respectively, are obtained through applying the isometric transform 
		$\mathcal{F}_n$,
		introduced in the proof of Lemma~\ref{lem:ka}, to $\left(\widetilde{k}^i\right)_{i=1}^n$ and 
		$\left(k^i\right)_{i=1}^n$, respectively. Thus, the estimate \eqref{eq:ktsa} also holds for the 
		vector-valued functions, i.e., 
		\begin{equation}
			\label{eq:ktsa2}
			\max_{(x,\xi)\in\mathcal{T}}\frac{1}{\sqrt{n}}\left\| \left(k^i(x,\xi)\right)_{i=1}^n - 
			\left(\widetilde{k}^i(x,\xi)\right)_{i=1}^n\right\|_{\mathbb{R}^n} \leq \delta + \varepsilon.
		\end{equation}
		In addition, from \eqref{eq:ka} in Lemma~\ref{lem:ka} we already have 
		\begin{equation}
			\label{eq:ktn+1a}
			\max_{(x,\xi)\in \mathcal{T}} \left|k^{n+1}(x,\xi) - \widetilde{k}^{n+1}(x,\xi)\right| \leq \delta.
		\end{equation}
		Now, setting $\Delta k^i = \widetilde{k}^i - k^i$ for $i=1,2,\ldots,n+1$, we have written 
		\eqref{eq:aUn}  as \eqref{eq:aU}, where the error term can be estimated using 
		\eqref{eq:ktsa2}, \eqref{eq:ktn+1a}, triangle inequality and Cauchy-Schwartz inequality as
		\begin{align}
			\label{eq:dke}
			\int\limits_0^1\left(\frac{1}{n}\sum_{i=1}^n\Delta k^i(1,\xi)u^i(t,\xi) +
			\Delta k^{n+1}(1,\xi)v(t,\xi)\right)d\xi & \leq \nonumber \\
			(2\delta + \varepsilon)\left\| \left(\begin{smallmatrix}
				\mathbf{u}(t,\cdot) \\ v(t,\cdot)
			\end{smallmatrix}\right)\right\|_E, &
		\end{align}
		where $\delta$ and $\varepsilon$ become arbitrarily small when $n$ is sufficiently large.
	\end{proof}
\end{lemma}

By Lemma~\ref{lem:aU}, the control law \eqref{eq:aUn} can be split into the part that 
exponentially 
stabilizes the large-scale $n+1$ system \eqref{eq:n+1}, \eqref{eq:nbcuy} and to the $\Delta$-part 
which we treat as a perturbation that 
becomes arbitrarily small when $n$ is sufficiently large. Thus, the stability of the $n+1$ system 
under the control law \eqref{eq:aU} can be established based on existing results for well-posed 
infinite-dimensional linear systems. The well-posedness of the $n+1$ system and a generic 
stability result for perturbed well-posed linear systems are presented in Propositions 
\ref{prop:wp} and~\ref{prop:wpstab}, respectively, in Appendix \ref{app:tech}.

\textit{Proof of Theorem~\ref{thm:stab}.} By Proposition~\ref{prop:wp} and \cite[Sect. 
10.1]{TucWeiBook}, we can translate the 
boundary-controlled PDE \eqref{eq:n+1}, \eqref{eq:nbcuy} into a well-posed abstract Cauchy 
problem $\dot{z}(t) = Az(t) + BU(t)$ on the Hilbert space $E$, 
where 
$z = (u^1,\ldots,u^n,v)^T$, $\dot{z}(t) = Az(t)$ corresponds to \eqref{eq:n+1} with the 
homogeneous boundary condition from  \eqref{eq:nbcuy} through the domain of $A$, and 
$BU(t)$ corresponds to the boundary control in \eqref{eq:nbcuy}.
Moreover, we introduce bounded linear operators $K$ and $\Delta K$ corresponding to taking 
inner products (on $E$) with $\left(k^i(1,\cdot)\right)_{i=1}^{n+1}$ and $\left(\Delta 
k^i(1,\cdot)\right)_{i=1}^{n+1}$ 
such that the 
control law \eqref{eq:aU} can be expressed as $U(t) = Kz(t) + \Delta Kz(t)$, 
where $\|\Delta K\|_{\mathcal{L}(E,\mathbb{R})} \leq 2\delta+\varepsilon$ becomes arbitrarily 
small when $n$ is sufficiently large by Lemma 
\ref{lem:aU}. Thus, the exponential stability of \eqref{eq:n+1}, \eqref{eq:nbcuy} under the control 
law \eqref{eq:aU} follows from \cite[Thm 3.2]{DiMVaz13} and Proposition~\ref{prop:wpstab}. As 
this control law is equivalent to \eqref{eq:aUn} by Lemma~\ref{lem:aU}, this concludes the proof 
of Theorem~\ref{thm:stab}.

\subsection{Proof of Theorem \ref{thm:stab} Using a Lyapunov Functional} \label{sec:lyap}

By following \cite[Lem. 3.1]{DiMVaz13}, we construct 
the Lyapunov functional with parameters $p,\delta_1 > 0$
\begin{equation}
	\label{eq:lyap}
	V(t) =\int\limits_{0}^{1} pe^{-\delta_1 x}\frac{1}{n}\sum_{i=1}^n
	\frac{\alpha^2_i(t,x)}{\lambda_i(x)}dx + \int\limits_0^1 \frac{1+x}{\mu(x)}\beta(t,x)^2dx,
\end{equation}
where $\alpha_i(t,x) = u_i(t,x)$ for $i=1,\ldots,n$ and 
\begin{align}
	\label{eq:beta}
	\beta(t,x) & = v(t,x) - \int\limits_0^x
	\frac{1}{n}\sum_{i=1}^nk^i(x,\xi)u^i(t,\xi)d\xi
	\nonumber \\
	& \qquad - \int\limits_0^x k^{n+1}(x,\xi)v(t,\xi)d\xi,
\end{align}
are the states of the target system after the exact backstepping
transformation. The dynamics of the target system is of the form
\cite[Sect. III.A]{DiMVaz13}
\begin{subequations}
	\label{eq:ts}%
	\begin{align}
	\alpha_t^i(t,x) +\lambda_i(x)\alpha_x^i(t,x) & = \nonumber \\
		\frac{1}{n}\sum_{j=1}^n	\sigma_{i,j}(x)\alpha^j(t,x)+W_i(x)\beta(t,x) \nonumber \\
		 +	\frac{1}{n}\sum_{j=1}^n
		\int\limits_0^xc_{i,j}(x,\xi)\alpha^j(t,\xi)d\xi
		+
		\int\limits_0^x\kappa_i(x,\xi)\beta(t,\xi)d\xi, & \\
		\beta_t(t,x)-\mu(x)\beta_x(t,x) & = 0,
	\end{align}
\end{subequations}
with boundary conditions $\alpha^i(t,0) = q_i\beta(t,0)$ 
for all $i=1,2,\ldots,n$, where $c_{i,j},\kappa_i$ are continuous on
$\mathcal{T}$\footnote{Under the exact backstepping controller, \eqref{eq:ts} would be 
	accompanied with boundary condition $\beta(t,1)=0$.}.
Employing the approximate 
control law \eqref{eq:aU}, the boundary condition for $\beta$ at $x=1$ is
\begin{align}
	\label{eq:atbc}
	\beta(t,1) & = \int\limits_0^1\frac{1}{n}\sum_{i=1}^n\Delta
	k^i(1,\xi)u^i(t,\xi)d\xi  \nonumber \\
	& \quad +  \int\limits_0^1\Delta k^{n+1}(1,\xi)v(t,\xi)d\xi \nonumber \\
	& = \left< \left(\begin{smallmatrix}
		\Delta\mathbf{k}(1,\cdot) \\ \Delta k^{n+1}(1,\cdot)
	\end{smallmatrix}\right), \left(\begin{smallmatrix}
		\mathbf{u}(t, \cdot) \\ v(t,\cdot)
	\end{smallmatrix}\right) \right>_E,
\end{align}
which also needs to be taken into account in the Lyapunov-based analysis.

Computing $\dot{V}(t)$ and integrating by
parts yields
\begin{equation}
	\label{eq:Vd}
	\resizebox{.99\columnwidth}{!}{$
	\begin{split}
	\dot{V}(t) & = \left[
	-pe^{-\delta_1x}\frac{1}{n}\|\pmb{\alpha}(t,x)\|^2_{\mathbb{R}^n}+(1+x)\beta(t,x)^2
	\right]_0^1\\ 
	& \quad -
	\int\limits_0^1 \left(\delta_1pe^{-\delta_1x}\frac{1}{n}\|\pmb{\alpha}(t,x)\|^2_{\mathbb{R}^n}
	+\beta(t,x)^2 \right)dx \\
	& \quad +
	\int\limits_0^1pe^{-\delta_1x}\frac{1}{n}\pmb{\alpha}(t,x)^T\pmb{\lambda}^{-1}(x)
	\frac{\pmb{\sigma}(x)}{n}\pmb{\alpha}(t,x)dx \\
	& \quad + 2 \int\limits_0^1
	\int\limits_0^xpe^{-\delta_1x}\frac{1}{n}\pmb{\alpha}(t,x)^T\pmb{\lambda}^{-1}(x)
	\frac{\mathbf{c}(x,\xi)}{n}\pmb{\alpha}(t,\xi)d\xi dx \\
	& \quad + 2
	\int\limits_0^1pe^{-\delta_1x}\frac{1}{n}\pmb{\alpha}(t,x)^T\pmb{\lambda}^{-1}(x)
	\mathbf{W}(x)\beta(t,x)dx \\
	& \quad + 2 \int\limits_0^1
	\int\limits_0^xpe^{-\delta_1x}\frac{1}{n}\pmb{\alpha}(t,x)^T
	\pmb{\lambda}^{-1}(x)\pmb{\kappa}(x,\xi) \beta(t,\xi)d\xi dx,
	\end{split}$}
\end{equation}
where we denote $\pmb{\alpha}=\left(\alpha_i\right)_{i=1}^n, 
\mathbf{W}=\left(W_i\right)_{i=1}^n,\pmb{\kappa} = \left(\kappa_i\right)_{i=1}^n$ and
$\pmb{\lambda}=\operatorname{diag}(\lambda_i)_{i=1}^n, \pmb{\sigma} = 
\left(\sigma_{i,j}\right)_{i,j=1}^n, \mathbf{c}=\left(c_{i,j}\right)_{i,j=1}^n$. Since all the
(individual components of the) parameters are continuous, they are
also uniformly bounded on compact sets, and hence, there exist some
$M_{\sigma}, M_W, M_c, M_{\kappa}>0$ such that 
\begin{subequations}
	\label{eq:lbd}%
	\begin{align}
		\max_{x\in[0,1]}\frac{\|\pmb{\sigma}(x)\|_{\infty}}{n}
		& \leq M_{\sigma}, \\
		\max_{x\in[0,1]}\|\mathbf{W}(x)\|_{\infty}
		& \leq M_W, \\
		\max_{(x,\xi)\in \mathcal{T}}
		\frac{\|\mathbf{c}(x,\xi)\|_{\infty}}{n}
		&\leq  M_c, \\
		\max_{(x,\xi)\in \mathcal{T}} \|\pmb{\kappa}(x,\xi)\|_{\infty}
		& \leq M_{\kappa}.
	\end{align}
\end{subequations}
Moreover, since $\pmb{\lambda}(x)$ is diagonal and uniformly bounded
away from zero by Assumption~\ref{ass:n+1}, there exists some
$M_{\lambda} > 0$ such that
$\max_{x\in[0,1]}\|\pmb{\lambda}^{-1}(x)\|_{\infty} \leq
M_{\lambda}$. Thus, we eventually get 
\begin{align}
	\label{eq:Vdest}
	\dot{V}(t) & \leq -\left(1 - p \frac{1}{n}\sum_{i=1}^nq_i^2
	\right)\beta(t,0)^2 + 2\beta(t,1)^2 \nonumber \\
	& \quad - \int\limits_0^1 \left(1 - pM_{\lambda}M_W -
	\frac{pM_{\lambda}M_{\kappa}}{\delta_1} \right)\beta(t,x)^2dx \nonumber \\
	& \quad - \int\limits_0^1 pe^{-\delta_1x}(\delta_1 -
	\widetilde{M})\frac{1}{n}\|\pmb{\alpha}(t,x)\|^2_{\mathbb{R}^n}dx,
\end{align}
where $\widetilde{M} = M_{\lambda}M_{\sigma} +
2\frac{M_{\lambda}M_c}{\delta_1} + M_{\lambda}M_W +
\frac{M_{\lambda}M_{\kappa}}{\delta_1}$\footnote{If $\beta(t,1)=0$, under the exact controller, 
then $\dot{V}(t)$ is negative definite for some sufficiently small $p>0$ and sufficiently large 
$\delta_1>0$.}. We have the following estimate 
	from
\eqref{eq:atbc} and \eqref{eq:dke}
\begin{equation}
	\label{eq:betaest}
	2\beta(t,1)^2 
	\leq 2(2\delta + \varepsilon)^2\left\| \left(\begin{smallmatrix}
		\mathbf{u}(t,\cdot) \\ v(t,\cdot)
	\end{smallmatrix}\right)\right\|_E^2,
\end{equation}
where $u^i = \alpha^i$ for all $i =1,2,\ldots,n$. Moreover, there exist
inverse kernels $\left( l^i \right)_{i=1}^{n+1}$ such that
\cite[Sect. III.A.3]{DiMVaz13} 
\begin{align}
	\label{eq:abtov}
	v(t,x) & = \beta(t,x) + \left< \left(\begin{smallmatrix}
		\mathbf{l}(x,\cdot) \\ l^{n+1}(x,\cdot)
	\end{smallmatrix}\right), \left(\begin{smallmatrix}
		\pmb{\alpha}(t, \cdot) \\ \beta(t,\cdot)
	\end{smallmatrix}\right) \right>_E,
\end{align}
where $\left( l^i \right)_{i=1}^{n+1}$ are continuous on $\mathcal{T}$,
and hence, also uniformly bounded. Thus, there exists some $M_l >
0$ such that 
\begin{equation}
	\label{eq:lest}
\max_{(x,\xi)\in \mathcal{T}}\max_{i\in \{1,2,\ldots,n+1\}}\left|l^i(x,\xi)\right| \leq M_l.
\end{equation}
By Jensen's inequality we thus have
\begin{align}
	\label{eq:uvtoab}
	\left\| \left(\begin{smallmatrix}
		\mathbf{u}(t,\cdot) \\ v(t,\cdot)
	\end{smallmatrix}\right)\right\|_E^2 & \leq 
	(1+M_{l})^2\left\| \left(\begin{smallmatrix}
		\pmb{\alpha}(t,\cdot) \\ \beta(t,\cdot)
	\end{smallmatrix}\right)\right\|_E^2,
\end{align}
which finally yields 
\begin{equation}
	\label{eq:b1est}
	2\beta(t,1)^2 \leq 2(2\delta+\varepsilon)^2(1+M_l)^2
	\left\| \left(\begin{smallmatrix}
		\pmb{\alpha}(t,\cdot) \\ \beta(t,\cdot)
	\end{smallmatrix}\right)\right\|_E^2.
\end{equation}
Relation \eqref{eq:b1est} implies that when $\delta$ and $\varepsilon$ are sufficiently small,
$\dot{V}(t)$ in \eqref{eq:Vdest} can be made negative definite, with a proper choice of 
$p,\delta_1$, despite the perturbation acting on $\beta(t,1)$.

To complete the proof, we note that $V(t)$ in \eqref{eq:lyap} corresponds to the weighted inner 
product $V(t) = \left\langle \left( 
\begin{smallmatrix}
	\pmb{\alpha}(t,\cdot) \\ \beta(t,\cdot)
\end{smallmatrix}
\right), \mathbf{P}(\cdot)\left( 
\begin{smallmatrix}
	\pmb{\alpha}(t,\cdot)\\ \beta(t,\cdot)
\end{smallmatrix}
\right) \right\rangle_E$, where $\mathbf{P}(x) = 
\operatorname{diag}\left(pe^{-\delta_1x}\pmb{\lambda}^{-1}(x), \frac{1+x}{\mu(x)}\right) > 0$ for 
all $x\in[0,1]$, by 
which there exists some $m_V,M_V > 0$ such that 
\begin{equation}
	\label{eq:Vest}
	m_V\left\| \left( 
	\begin{smallmatrix}
		\pmb{\alpha}(t,\cdot) \\ \beta(t,\cdot)
	\end{smallmatrix}
	\right) \right\|^2_E \leq V(t) \leq M_V\left\| \left( 
	\begin{smallmatrix}
		\pmb{\alpha}(t,\cdot) \\ \beta(t,\cdot)
	\end{smallmatrix}
	\right) \right\|^2_E.
\end{equation}
Moreover, from \eqref{eq:Vdest} and \eqref{eq:b1est} we get $\dot{V}(t) \leq 
-c_V\left\| \left( 
\begin{smallmatrix}
	\pmb{\alpha}(t,\cdot) \\ \beta(t,\cdot)
\end{smallmatrix}
\right) \right\|^2_E$ with $c_V = \min\left\{1-pM_\lambda M_W - \frac{pM_\lambda 
	M_\kappa}{\delta_1}, 
pe^{-\delta_1}\left(\delta_1-\widetilde{M}\right)\right\}-2(2\delta+\varepsilon)^2(1+M_l^2)$, which 
can be made positive for a large $\delta_1$ and small $p,\delta,\varepsilon$.
Thus, we have 
\begin{equation}
	\label{eq:abdecay}
	\left\| \left( 
	\begin{smallmatrix}
		\pmb{\alpha}(t,\cdot) \\ \beta(t,\cdot)
	\end{smallmatrix}
	\right) \right\|^2_E \leq \frac{m_V}{M_V}\left\| \left( 
	\begin{smallmatrix}
		\pmb{\alpha}(0,\cdot) \\ \beta(0,\cdot)
	\end{smallmatrix}
	\right) \right\|^2_E\exp\left(-\frac{c_V}{M_V}t\right).
\end{equation}
Finally, using \eqref{eq:uvtoab} and an analogous estimate for the forward transform
\eqref{eq:beta} as
$\left\| \left( 
\begin{smallmatrix}
	\pmb{\alpha}\\ \beta
\end{smallmatrix}
\right) \right\|_E \leq  (1+M_k)\left\| \left( 
\begin{smallmatrix}
	\pmb{u} \\ v
\end{smallmatrix}
\right) \right\|_E 
$,
for some $M_k>0$ such that $\displaystyle \max_{(x,\xi)\in\mathcal{T}}\max_{i\in 
\{1,2,\ldots,n+1\}}\left|k^i(x,\xi)\right|\leq M_k$, we obtain
\begin{equation}
	\label{eq:uvdecay}
	\left\| \left( 
	\begin{smallmatrix}
		\mathbf{u}(t,\cdot) \\ v(t,\cdot)
	\end{smallmatrix}
	\right) \right\|^2_E \leq \widetilde{M}_V\left\| \left( 
	\begin{smallmatrix}
		\mathbf{u}(0,\cdot) \\ v(0,\cdot)
	\end{smallmatrix}
	\right) \right\|^2_E\exp\left(-\frac{c_V}{M_V}t\right),
\end{equation}
where $\widetilde{M}_V = (1+M_l)^2\frac{m_V}{M_V}(1+M_k)^2$, which completes the proof.

\begin{remark}
	Technically, the Lyapunov-based arguments presented apply to classical solutions, the 
	existence of which 
	can be guaranteed for any initial conditions $u_0^i, v_0 \in H^1([0,1];\mathbb{R})$ that satisfy 
	the compatibility conditions $u_0^i(0) = q_iv_0(0)$ and $v_0(1) = U(0)$ \cite[Prop. 
	10.1.8]{TucWeiBook}. However, as noted in 
	\cite[Sect. 2.1.3]{BasCorBook}, for any weak solution there exists a sequence of classical 
	solutions which converges to the weak solution in $E$, and hence, the decay 
	estimate 
	\eqref{eq:uvdecay} also applies to weak solutions.
\end{remark}

\section{Numerical Example} \label{sec:sim}

\subsection{Stabilization via Approximate Kernels} \label{sim:a}

As an example, consider an $n+1$ system
\eqref{eq:n+1}, \eqref{eq:nbcuy} with parameters 
$\mu(x) = 1$ and
\begin{subequations}
	\label{eq:exp}%
	\begin{align}
		\lambda_i(x)
		& = 1, \\
		\sigma_{i,j}(x)
		& = x^3(x+1)\left(\frac{i}{n}-\frac{1}{2} \right)
		\left( \frac{j}{n} - \frac{1}{2} \right), \\
		W_i(x)
		& = x(x+1)e^x\left( \frac{i}{n} - \frac{1}{2} \right), \\
		\theta_i(x)
		& = -70e^{x \frac{35}{\pi^2}}\frac{i}{n}\left(
		\frac{i}{n}-1 \right), \\
		q_i & = \cos\left(2\pi \frac{i}{n}\right),
	\end{align}%
\end{subequations}
for $i,j = 1,\ldots,n$ such that continuous functions satisfying \eqref{eq:afn} can be constructed as 
\begin{subequations}
	\label{eq:expc}%
	\begin{align}
		\lambda(x,y)
		& = 1, \\
		\sigma(x,y,\eta)
		& = x^3(x+1)\left( y-\frac{1}{2} \right)\left( \eta-\frac{1}{2}
		\right), \\
		W(x,y)
		& = x(x+1)e^x \left( y - \frac{1}{2} \right), \\
		\theta(x,y)
		& = -70e^{x \frac{35}{\pi^2}} y \left( y - 1 \right), \\
		q(y) & = \cos(2\pi y).
	\end{align}
\end{subequations}
The latter parameter values correspond to the example considered in
\cite[Sect. VII]{AllKrs23arxiv}, where it is shown that the solutions
of the corresponding continuum kernel equations are given by 
\begin{subequations}
	\label{eq:exkc}%
	\begin{align}
		k(x,\xi,y) & = 35y(y-1)e^{2\xi \bar{k}(x,\xi)}, \label{eq:exkca} \\
		\bar{k}(x,\xi) & = \frac{35}{2\pi^2}.
	\end{align}
\end{subequations}
We note that while a closed-form solution exists to the continuum kernel
equations \eqref{eq:kc}, \eqref{eq:kcbc} with parameters \eqref{eq:expc}, we were not able to 
solve the corresponding $n+1$ kernel
equations \eqref{eq:kn}, \eqref{eq:knbc} with parameters \eqref{eq:exp} in closed-form when
$n\in \mathbb{N}$ is arbitrary, nor do we expect that a closed-form
solution can be constructed (even for small $n$). This is also consistent with, e.g., 
\cite{VazKrs14}, in which 
an explicit solution is possible to obtain for the specific case $n=1$ and for spatially invariant 
parameters \eqref{eq:exp}. Regardless, in this particular example, the
continuum approximation significantly simplifies the computation of the stabilizing control kernels.

We simulate the $n+1$ system \eqref{eq:n+1}, \eqref{eq:nbcuy} with parameters \eqref{eq:exp} 
under the control law \eqref{eq:aUn} computed based on the continuum kernels \eqref{eq:exkc}. 
Various values of $n$ are considered to illustrate the behavior of the closed-loop system as $n$ 
increases. In fact, the closed-loop system is stable for any $n\geq 2$, but the performance 
is improved for larger $n$. However, when $n=1$, the system \eqref{eq:n+1}, \eqref{eq:nbcuy} is 
open-loop stable and $k^1(x,\xi) = k^2(x,\xi) = 0$ is the solution to the kernel 
equations \eqref{eq:kn}, \eqref{eq:knbc}, in which case the approximate control law 
\eqref{eq:aUn} destabilizes the system 
(because $\bar{k}\neq 0$, in contrast to $k^2$).

\begin{figure}[!t]
	\includegraphics[width=\columnwidth]{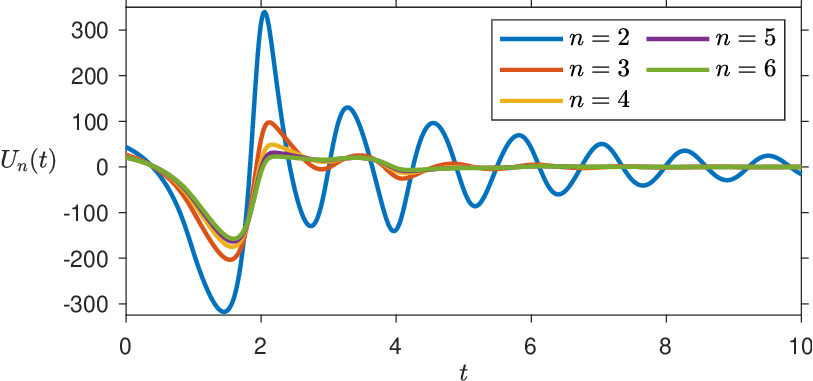}
	\caption{The controls $U(t)$ based on the approximate control law \eqref{eq:aUn} when 
		$n=2,\ldots,6$.}
	\label{fig:U1}
\end{figure}

\begin{figure}[htbp]
	\includegraphics[width=\columnwidth]{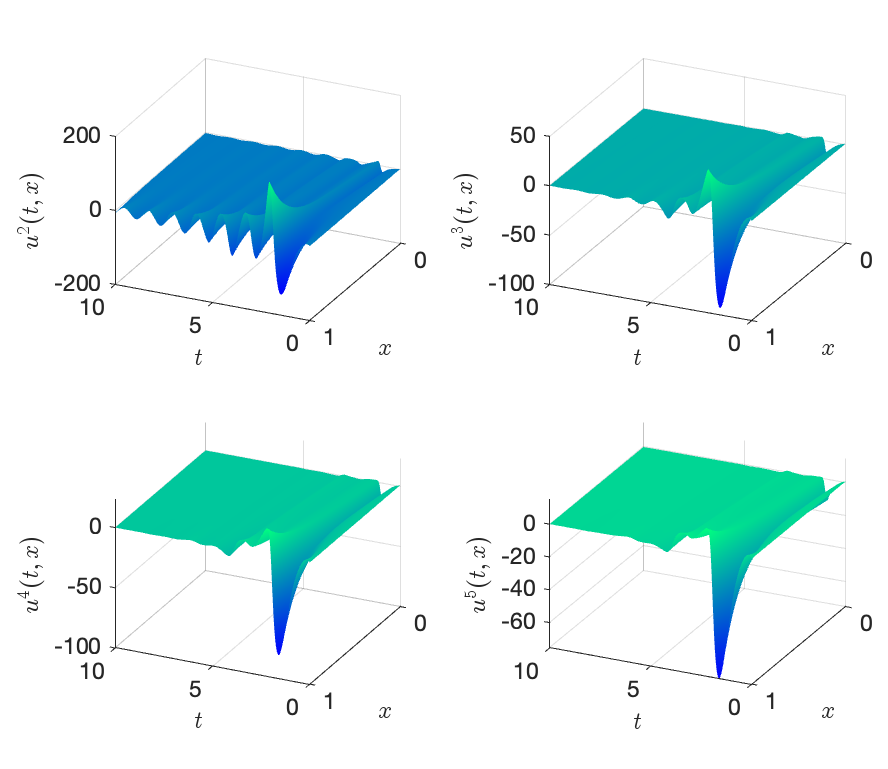}
	\caption{The solution component $u^n(t,x)$ when	$n=2,\ldots,5$.}
	\label{fig:sol1}
\end{figure}

In the simulations, the system \eqref{eq:n+1}, \eqref{eq:nbcuy} is approximated using finite 
differences with $256$ grid points in $x\in[0,1]$. The ODE resulting from the finite-difference 
approximation is solved using \verb+ode45+ in MATLAB. The initial conditions for all $n$ are 
$u_0^i(x) = q_i$ for $i=1,\ldots,n$ and $v_0(x) = 1$, for all $x\in[0,1]$. Fig.~\ref{fig:U1} displays 
the control law 
\eqref{eq:aUn} when $n=2,\ldots,6$. We note that the control law also acts as a weighted average 
of the solution to \eqref{eq:n+1}, \eqref{eq:nbcuy}, i.e., we can also assess the exponential decay 
rate of the solutions based on $U(t)$. However, as $k(x,\xi,1) = 0$ in \eqref{eq:exkc}, the 
component $u^n(t,x)$ does not affect the control in any way. Therefore, $u^n(t,x)$ is displayed 
separately in Fig.~\ref{fig:sol1} for $n=2,\ldots,5$, which shows that also $u^n(t,x)$ decays to 
zero as $t\to\infty$. 

\begin{figure}[!t]
	\includegraphics[width=\columnwidth]{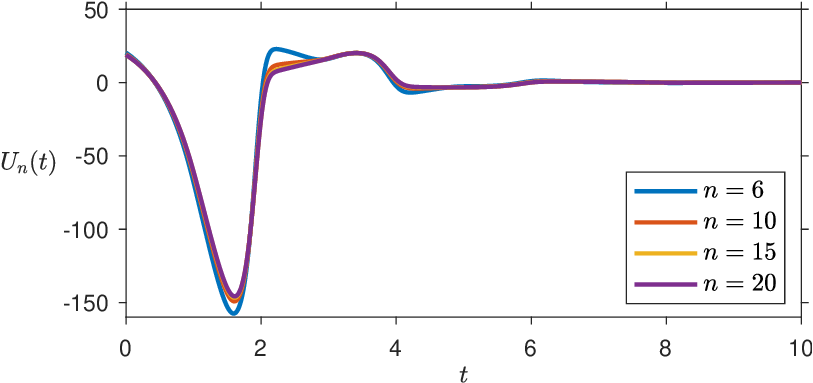}
	\caption{The controls $U(t)$ based on the approximate control law \eqref{eq:aUn} 
		when $n=6,10,15,20$.}
	\label{fig:U2}
\end{figure}

\begin{figure}[!t]
	\includegraphics[width=\columnwidth]{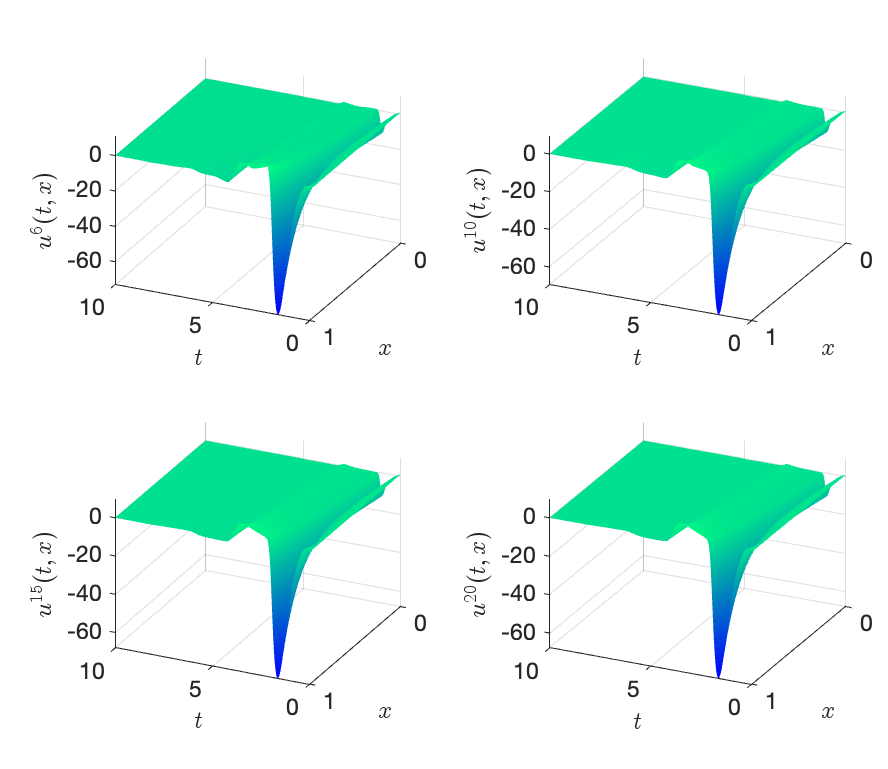}
	\caption{The solution component $u^n(t,x)$ when $n=6,10,15,20$.}
	\label{fig:sol2}
\end{figure}

Fig.~\ref{fig:U1} and Fig.~\ref{fig:sol1} show that the approximate control law based on the 
continuum kernels \eqref{eq:exkc} is indeed stabilizing already when $n=2$, even if the rate of 
decay is very slow. However, Fig.~\ref{fig:U1} and Fig.~\ref{fig:sol1} show that the closed-loop 
performance significantly 
improves when $n$ becomes larger, and in Fig.~\ref{fig:U1} the controls for $n=5$ and $n=6$ are 
virtually indistinguishable. However, as we consider larger values of $n$ separately in Fig.
\ref{fig:U2} and Fig.~\ref{fig:sol2}, some changes in $U(t)$ and $u^n(t,\xi)$ are still noticeable 
between 
$n=6$ and $n=10$. Regardless, in all studied cases beyond $n>5$, the controls, along with the 
solutions, have practically converged to zero by $t=6$. We note that it takes the control input 
$1/\mu + 1/\lambda = 2$ time units to traverse through the system \eqref{eq:n+1}, 
\eqref{eq:nbcuy} for any $n$, which is why, e.g., the state component $u^n(t,\xi)$, may grow 
rapidly in the beginning of the simulation, as seen in Fig.~\ref{fig:sol1} and Fig.~\ref{fig:sol2}, 
before 
getting stabilized by the controls. However, due to the in-domain coupling between $u$ and $v$ 
in \eqref{eq:n+1}, the controls do affect the $u^i$ components through $v$ already before 
entering the $u^i$ channels at time $t=1/\mu = 1$.

Overall, the simulations demonstrate that the approximate control law \eqref{eq:aUn} based on 
the continuum kernels \eqref{eq:exkc} exponentially stabilizes the $n+1$ system \eqref{eq:n+1}, 
\eqref{eq:nbcuy} when $n$ is sufficiently large, and that the approximation error of the control law 
decreases as $n$ increases. Thus, the simulations are well in accordance with the theoretical 
results. Moreover, in this example, the approximate control law has good performance already for 
very moderate $n$, showing that the sufficiently large $n$ appearing in the theoretical results 
may be, in practice, relatively small. However, one should not expect this to be the case in 
general, as this is 
dependent on the parameters of both the system \eqref{eq:n+1}, \eqref{eq:nbcuy} and its 
continuum approximation \eqref{eq:inf}, \eqref{eq:cbcuy}.

\subsection{Comparison with Exact Control Kernels} \label{sim:comp}

In this subsection, we compare the approximate control gain kernels with the exact kernels, as 
well as we compare performance of the closed-loop systems under the approximate and exact 
kernels. However, as 
we are not aware of the existence of closed-form solutions to \eqref{eq:kn}, \eqref{eq:knbc}, we 
have to find the solution implicitly, and hence, the presented comparisons are between the 
closed-form continuum kernels \eqref{eq:exkc} and numerical approximations of the exact kernels 
obtained from \eqref{eq:kn}, \eqref{eq:knbc}. Regardless, we refer to the numerical solution to 
\eqref{eq:kn}, \eqref{eq:knbc} as the exact kernels.

\begin{figure}[!t]
	\includegraphics[width=\columnwidth]{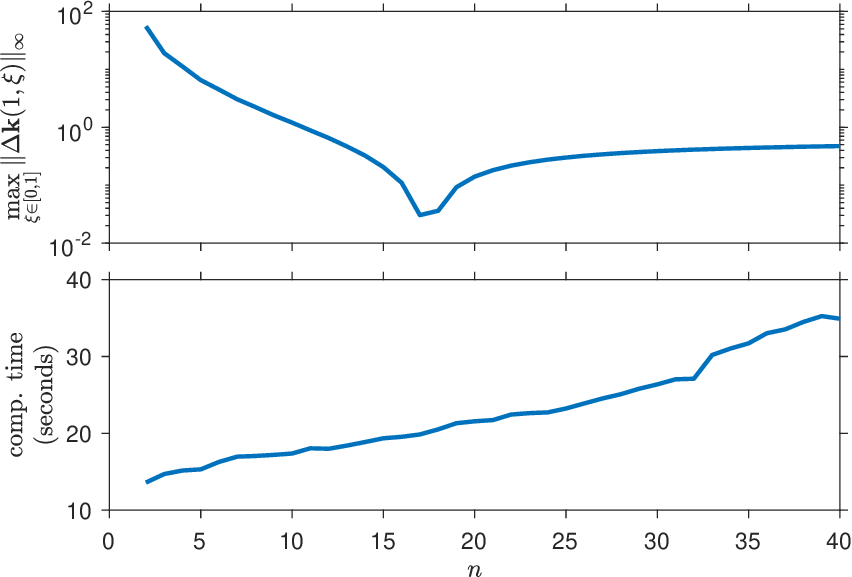}
	\caption{Top: Approximation error of the exact control kernels by the continuum control 
	kernels. Bottom: Computational time to solving the exact control kernels PDEs.}
	\label{fig:kercomp}
\end{figure}

In Fig.~\ref{fig:kercomp}, the evolution of $\displaystyle \max_{\xi\in[0,1]}\|\Delta 
\mathbf{k}(1,\xi)\|_\infty$ is shown for $n=2,3,\ldots,40$, where we use the boldface notation for 
the vector $\left(\Delta k^i\right)_{i=1}^{n+1}$. It can be seen that the error norm decays (even 
exponentially) as $n$ increases, until $n=17$ after which the accuracy of our numerical procedure 
to 
solving the exact kernel equations \eqref{eq:kn}, \eqref{eq:knbc} starts to deteriorate. That is, the 
increase in the error norm around $n=20$ is due to numerical inaccuracies in solving the exact 
kernel equations \eqref{eq:kn}, \eqref{eq:knbc}, as in theory the error norm should tend towards 
zero as $n$ increases. Fig.~\ref{fig:kercomp} additionally shows the computational time required 
(by our numerical procedure) to solve \eqref{eq:kn}, \eqref{eq:knbc}, which appears to grow 
linearly with respect to $n$. However, solving \eqref{eq:kn}, \eqref{eq:knbc} additionally requires 
spatial discretization of the domain $\mathcal{T}$, where we simply use finite differences with 
257 grid points in both $x$ and $\xi$. We expect that \eqref{eq:kn}, \eqref{eq:knbc} can be solved 
more accurately using more sophisticated methods, such as finite element methods (FEM), and a 
finer discretization of 
$\mathcal{T}$, but this will likely result in longer computational times as well, since the complexity 
of FEM algorithms is proportional to the number of nodes in the finite element mesh 
\cite{FarShmCADSM11}. We note here that the computational time of computing the approximate 
kernels is virtually invariant to $n$, as such computations rely only on evaluation of 
\eqref{eq:exkca} at the given points in $y$.

\begin{figure}[htbp]
	\includegraphics[width=\columnwidth]{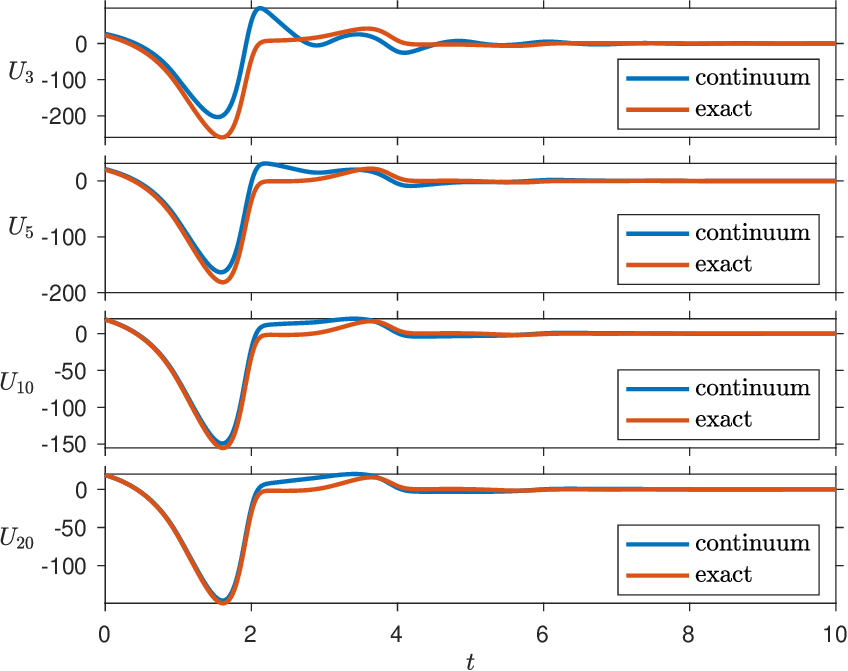}
	\caption{Comparison of the control efforts with the continuum (approximate) kernels and the 
		exact control kernels for $n=3,5,10,20$.}
	\label{fig:solcomp}
\end{figure}

In Fig.~\ref{fig:solcomp}, we compare the control efforts obtained by using the approximate 
continuum 
kernels and the exact kernels obtained from  \eqref{eq:kn}, \eqref{eq:knbc} for $n=3,5,10,20$. It 
can be seen that qualitatively the controls computed based on the exact kernels seem to behave 
the same 
despite of $n$. Moreover, as expected from 
Fig.~\ref{fig:U1} and Fig.~\ref{fig:U2}, the approximate control laws computed based on the 
continuum 
kernels \eqref{eq:exkc} get closer to the exact controls as $n$ increases. However, as already 
seen in Fig.~\ref{fig:kercomp}, we can numerically demonstrate this convergence only up to some 
finite accuracy. Regardless, it can be seen from Fig.~\ref{fig:solcomp} that the differences in 
closed-loop performance are small already at $n=10$.

\section{Approximation of the Solutions to the Large-Scale System by a Continuum} 
\label{sec:part}

In this section, we show that the continuum system \eqref{eq:inf}, \eqref{eq:cbcuy} acts as an 
approximation of the $n+1$ system \eqref{eq:n+1}, \eqref{eq:nbcuy}, when $n$ is sufficiently 
large. In particular, a sequence of solutions to \eqref{eq:n+1}, \eqref{eq:nbcuy} can approximate 
the solution to \eqref{eq:inf}, \eqref{eq:cbcuy}, provided that the sequence of data of 
\eqref{eq:n+1}, \eqref{eq:nbcuy} approximates the data of \eqref{eq:inf}, \eqref{eq:cbcuy} 
(including initial conditions) to arbitrary accuracy.

\begin{theorem}
	\label{thm:ua}
	Consider an $n+1$ system \eqref{eq:n+1}, \eqref{eq:nbcuy}, with parameters $\lambda_i, \mu,
	W_i, \theta_i, \sigma_{i,j}, q_i$ for $i,j=1,2,\ldots,n$ satisfying 
	Assumption~\ref{ass:n+1}, initial conditions $v_0, u_0^i \in L^2([0,1];\mathbb{R})$ for 
	$i=1,2,\ldots,n$, and input $U \in L^2_{\rm{loc}}([0,+\infty); \mathbb{R})$. Construct a 
	continuum 
	system \eqref{eq:inf}, \eqref{eq:cbcuy}, with parameters $\lambda,\mu,W, 
	\theta, 	\sigma,q$ that satisfy Assumption~\ref{ass:inf} and \eqref{eq:afn}, and equip 
	\eqref{eq:inf}, \eqref{eq:cbcuy} with initial conditions $u_0,v_0$ and input $U$, such that 
	$u_0$ is continuous in $y$ and satisfies
	\begin{equation}
		\label{eq:u0c}
		u_0(x,i/n) = u_0^i(x),
	\end{equation}
	for $i=1,2,\ldots,n$\footnote{Such functions can always be constructed as, e.g.,
		per Footnote \ref{fn:afn}.}. Sample the solution 
	$(u,v)$ to the resulting PDE system \eqref{eq:inf}, 
	\eqref{eq:cbcuy} for these data, into a vector-valued function 
	$\left((\widetilde{u_i})_{i=1}^n,\widetilde{v}\right)$ as $\widetilde{\mathbf{u}}(t,x) = 
	\mathcal{F}_n^*u(t,x,\cdot)$\footnote{The transform 
		$\mathcal{F}_n$ is introduced in the proof of Lemma~\ref{lem:se}, and its adjoint 
		$\mathcal{F}_n^*$ satisfies \eqref{eq:Fns}.}  and $\widetilde{v}(t,x) = 
	v(t,x)$, pointwise for all $t \geq 0$ and almost all $x\in[0,1]$. On any compact interval 
	$t\in[0,T]$, for any given $T > 0$, we have 
	\begin{equation}
		\label{eq:unest}
		\max_{t\in[0, T]}\left\| \left( 
		\begin{smallmatrix}
			\mathbf{u}(t) \\ v(t)
		\end{smallmatrix}
		\right) -  \left( 
		\begin{smallmatrix}
			\widetilde{\mathbf{u}}(t) \\ \widetilde{v}(t)
		\end{smallmatrix}
		\right) \right\|_E \leq \delta_2 + \varepsilon_2 + \delta_3 + \varepsilon_3,
	\end{equation}
	where $\varepsilon_2,\varepsilon_3, \delta_2,\delta_3> 0$ become arbitrarily small when $n$ is 
	sufficiently large.
	
	\begin{proof}
		Firstly, we note that \eqref{eq:inf}, \eqref{eq:cbcuy} is well-posed on $E_c = 
		L^2([0,1];L^2([0,1];\mathbb{R}))\times L^2([0,1];\mathbb{R})$ under 
		Assumption~\ref{ass:inf} as
		shown in Proposition~\ref{prop:cwp} in Appendix~\ref{app:wp}. 
		That is, for any initial conditions $u_0 \in L^2([0, 1];L^2([0,1];\mathbb{R})), v_0 \in 
		L^2([0,1];\mathbb{R})$ and input $U \in L_{\rm{loc}}^2([0, +\infty);\mathbb{R})$, the solution 
		$(u,v)$ is well-defined, unique, continuous in time, and depends continuously on the 
		data of the problem on any compact interval $t \in [0, T]$, for any given $T > 0$. In 
		particular, 
		there exist families of bounded linear operators $\mathbb{T}_t, \Phi_t$ for $t\geq 0$, 
		depending continuously on the parameters $\lambda, \mu, \sigma, \theta, W, q$, such that 
		the 
		solution to  \eqref{eq:inf}, \eqref{eq:cbcuy} is given by \cite[Prop. 4.2.5]{TucWeiBook}
		\begin{equation}
			\label{eq:ucsol}
			\left(\begin{smallmatrix}
				u(t) \\ v(t)
			\end{smallmatrix}\right) = \mathbb{T}_t	\left(\begin{smallmatrix}
				u_0 \\ v_0
			\end{smallmatrix}\right) + \Phi_tU,
		\end{equation}
		which satisfies $\left(\begin{smallmatrix}
			u \\ v 	\end{smallmatrix}\right) \in C([0,+\infty); E_c)$. The well-posedness and regularity 
			of 
		solutions to \eqref{eq:n+1}, \eqref{eq:nbcuy} have already been established in 
		Proposition~\ref{prop:wp} in Appendix~\ref{app:tech} and Remark~\ref{rem:n+1}.
		
		Secondly, based on the parameters $\lambda_i, W_i, \theta_i, \sigma_{i,j}, q_i$ for 
		$i,j=1,2,\ldots,n$ of 
		\eqref{eq:n+1}, \eqref{eq:nbcuy}, we can construct the step functions $\lambda^n, \sigma^n, 
		\theta^n, W^n,  q^n$ as in \eqref{eq:sfn} such that they satisfy \eqref{eq:sa} for any 
		$\varepsilon_1 > 0$ when $n$ is sufficiently large. Moreover, we can transform the solution to
		\eqref{eq:n+1}, \eqref{eq:nbcuy} into a step function 
		(in $y$) by applying the transform $\mathcal{F}$ introduced in the proof of 
		Lemma~\ref{lem:se}, i.e., 
		\begin{equation}
			\label{eq:undef}
			\left( \begin{smallmatrix}
				u^n(t,x,\cdot) \\ v^n(t,x)
			\end{smallmatrix}\right) = \mathcal{F} \left(\begin{smallmatrix}
				\mathbf{u}(t,x) \\ v(t,x)
			\end{smallmatrix}\right).
		\end{equation}
		That is, $u^n$ is defined piecewise in $y$, for $i=1,2,\ldots,n$, as 
		$u^n(t,x,y) = u^i(t,x)$, $y\in ((i-1)/n,i/n]$,
		pointwise for all $t\in[0,T]$ and almost all $x\in[0,1]$. Then $(u^n,v^n)$ is the solution to 
		\eqref{eq:inf}, \eqref{eq:cbcuy} for parameters $\lambda^n, \mu, \sigma^n, 
		\theta^n, W^n,  q^n$, input $U$, and initial conditions $v_0, u_0^n(x,\cdot) = 
		\mathcal{F}_n\mathbf{u}_0(x)$. 
		This can be verified by applying the transform $\mathcal{F}$ to \eqref{eq:n+1}, 
		\eqref{eq:nbcuy} 
		(from the left), for each $t\geq 0$ and almost all $x\in[0,1]$, which results in $(u^n,v^n)$ 
		satisfying 
		\eqref{eq:inf}, \eqref{eq:cbcuy} for the stated parameters, initial conditions, and input. (This 
		follows in a similar manner with the respective part of the proof of Lemma~\ref{lem:se}.) As 
		\eqref{eq:n+1}, \eqref{eq:nbcuy} is well-posed, the transformed solution \eqref{eq:n+1}, 
		\eqref{eq:nbcuy} is the well-posed solution to the transformed (by $\mathcal{F}$) PDE 
		system. 
		In particular, there 
		exist families of operators $\mathbb{T}_t^n, \Phi_t^n$ for $t\geq 0$, depending continuously 
		on 
		the parameters $\lambda^n, \mu, \sigma^n, \theta^n, W^n, q^n$, such that the solution to 
		\eqref{eq:inf}, \eqref{eq:cbcuy} for these parameters is given by
		\begin{equation}
			\label{eq:ucnsol}
			\left(\begin{smallmatrix}
				u^n(t) \\ v^n(t)
			\end{smallmatrix}\right) = \mathbb{T}^n_t	\left(\begin{smallmatrix}
				u_0^n \\ v_0
			\end{smallmatrix}\right) + \Phi^n_tU,
		\end{equation}
		which satisfies $\left(\begin{smallmatrix}
			u^n\\ v^n \end{smallmatrix}\right) \in C([0,+\infty); E_c)$.
		
		Thirdly, consider the difference of the solutions (to  
		\eqref{eq:inf}, \eqref{eq:cbcuy}) \eqref{eq:ucsol} and \eqref{eq:ucnsol},   under the 
		respective parameters and initial conditions, 
		and the same inputs on the interval $t \in [0, T]$. We get, for each $t\in[0,T]$,
		\begin{align}
			\label{eq:ucest}
			\left\|  \left(\begin{smallmatrix}
				u(t) \\ v(t)
			\end{smallmatrix}\right)  -  \left(\begin{smallmatrix}
				u^n(t) \\ v^n(t)
			\end{smallmatrix}\right) \right\|_{E_c} & \leq \left\|\left( \mathbb{T}_t - 
			\mathbb{T}_t^n\right)\left(\begin{smallmatrix}
				u_0^n \\ v_0 	\end{smallmatrix}\right)\right\|_{E_c} \nonumber \\
			& \quad + \left\|\mathbb{T}_t\right\|_{\mathcal{L}(E_c)}\left\|\left(\begin{smallmatrix}
				u_0 \\ v_0 	\end{smallmatrix}\right) - \left(\begin{smallmatrix}
				u_0^n \\ v_0 	\end{smallmatrix}
			\right)\right\|_{E_c} \nonumber \\
			& \quad + \left\|\left(\Phi_t - \Phi_t^n\right)U\right\|_{E_c} \nonumber \\
			& \leq \delta_2 + \varepsilon_2 + \delta_3,
		\end{align}
		where, for any given $T, u_0,v_0$, and $U$, the constants $\delta_2,\varepsilon_2,\delta_3$ 
		can 
		be made 
		arbitrarily small by taking $n$ sufficiently large. That is, $\delta_2$ and 
		$\delta_3$\footnote{Note 
			that $\delta_2$ and $\delta_3$ are the maxima over $[0,T]$ of the respective constants 
			satisfying \eqref{eq:ucest} for each $t\in[0,T]$.} become small 
		due to the solution 
		to \eqref{eq:inf}, \eqref{eq:cbcuy} depending continuously on the parameters, initial 
		conditions, 
		and input of the problem, i.e., $\delta_2$ and $\delta_3$ depend continuously on 
		$\varepsilon_1$ 
		in \eqref{eq:sa} such that $\delta_2, \delta_3 \to 0$ as $\varepsilon_1\to 0$. The second term 
		$\varepsilon_2$ becomes small as  
		$\left\|\mathbb{T}_t\right\|_{\mathcal{L}(E_c)}$ is uniformly bounded on $t\in[0,T]$ by 
		Proposition~\ref{prop:cwp}; while the 
		difference $\|u_0-u_0^n\|_{L^2([0,1]^2;\mathbb{R})}$ can be made arbitrarily small due to 
		\eqref{eq:u0c} and 
		$u_0^n(x,\cdot) = \mathcal{F}_n\mathbf{u}_0(x)$, analogously to \eqref{eq:sa}.
		
		Finally, set $ \widetilde{\mathbf{u}}(t,x) = \mathcal{F}_n^*u(t,x,\cdot)$ and $\widetilde{v}=v$ 
		for 
		all $t \in [0,T]$ and almost all $x\in[0,1]$, i.e., 
		each component of  $\widetilde{\mathbf{u}}$ is 
		\begin{equation}
			\label{eq:utildedef}
			\widetilde{u}^i	(t,x)= n\int\limits_{(i-1)/n}^{i/n}u(t,x,y)dy,
		\end{equation}
		which is the mean value of $u(t,x,\cdot)$ over an interval of length $1/n$ in $y$. Since 
		$\mathcal{F}$ is 
		an isometry, we have, for each $t\in[0,T]$,
		\begin{align}
			\label{eq:usest}
			\left\| \left(\begin{smallmatrix}
				\mathbf{u}(t) \\ v(t)
			\end{smallmatrix}\right) - \left(\begin{smallmatrix}
				\widetilde{\mathbf{u}}(t) \\ \widetilde{v}(t)
			\end{smallmatrix}\right) \right\|_E 
			& =  \left\| \mathcal{F}\left( \left(\begin{smallmatrix}
				\mathbf{u}(t) \\ v(t)
			\end{smallmatrix}\right) - \left(\begin{smallmatrix}
				\widetilde{\mathbf{u}}(t)\\ \widetilde{v}(t)
			\end{smallmatrix}\right) \right)\right\|_{E_c} \nonumber \\
			& \leq  \left\|\left(\begin{smallmatrix}
				u^n(t) \\ v^n(t)
			\end{smallmatrix}\right)  -   \left(\begin{smallmatrix}
				u(t) \\ v(t)
			\end{smallmatrix}\right)  \right\|_{E_c} \nonumber \\
			& \quad +  \left\|  \left(\begin{smallmatrix}
				u(t) \\ v(t)
			\end{smallmatrix}\right)  -  \mathcal{F}\mathcal{F}^*  \left(\begin{smallmatrix}
				u(t) \\ v(t)
			\end{smallmatrix}\right)\right\|_{E_c} \nonumber \\
			& \leq \delta_2 + \varepsilon_2 + \delta_3 + \varepsilon_3,
		\end{align}
		where we used the triangle inequality (added and subtracted $(u,v)$), definitions 
		\eqref{eq:undef} 
		and \eqref{eq:utildedef}, and \eqref{eq:ucest}. 
		Moreover, $\varepsilon_3 > 0$ can be made arbitrarily small, as the solution $(u,v)$ is 
		uniformly 
		bounded on $t\in[0,T]$ and the step function $ \mathcal{F}\mathcal{F}^*  
		\left(\begin{smallmatrix}
			u(t) \\ v(t) \end{smallmatrix}\right)$ can approximate $\left(\begin{smallmatrix}
			u(t) \\ v(t) \end{smallmatrix}\right)$ to arbitrary accuracy (for almost all $x,y\in[0,1]$, 
			uniformly 
		in $t\in[0,T]$) when $n$ is 
		sufficiently large, i.e., when the interval length $1/n$ becomes sufficiently small (by 
		mean-value approximation, see, e.g., \cite[Sect. 1.6]{TaoBook11}). Thus, the 
		maximum of \eqref{eq:usest} over $t\in[0,T]$
		can be made arbitrarily small by taking $n$ sufficiently large, which yields \eqref{eq:unest}. 
	\end{proof}
\end{theorem}

\begin{remark}
	The conclusion of Theorem~\ref{thm:ua} remains valid even if the input $U$ was not exactly the 
	same for both PDE systems. That is, rewriting \eqref{eq:ucest} with $U$ replaced by some 
	$\bar{U}$ in \eqref{eq:ucsol} results in an additional error term 
	$\|\Phi_t\|_{\mathcal{L}(L^2([0, T];\mathbb{R}); 
		E_c)}\|U-\bar{U}\|_{L^2([0,T],\mathbb{R})}$, which can be made arbitrarily small by assuming 
	$\|U-\bar{U}\|_{L^2([0,T],\mathbb{R})}$ is sufficiently small, since the operator norm of 
	$\Phi_t$ is uniformly bounded on compact intervals $t\in[0,T]$ (based on 
	Proposition~\ref{prop:cwp} in Appendix~\ref{app:wp}). 
\end{remark}

\section{The Limiting Case $n\to\infty$} \label{sec:ntoinf}

While the system \eqref{eq:n+1}, \eqref{eq:nbcuy} and the kernel equations 
\eqref{eq:kn}, \eqref{eq:knbc} are well-defined for any finite $n$, we can consider the limiting 
case 
$n\to\infty$ through interpreting these as the respective continuum limits of system \eqref{eq:inf}, 
\eqref{eq:cbcuy} and kernel equations \eqref{eq:kc}, \eqref{eq:kcbc}. This can be formally proved 
as follows, utilizing the results of Sections \ref{sec:stab} and \ref{sec:part}. As regards the case of 
the kernels, it follows from \eqref{eq:ka} that 
	\begin{subequations}
	\label{eq:kainf}
	\begin{align}
		\lim_{n\to\infty}\max_{(x,\xi)\in \mathcal{T}}\|k(x,\xi,\cdot) - 
		k^n(x,\xi,\cdot)\|_{L^2([0,1];\mathbb{R})} & = 0, \\
		\lim_{n\to\infty}\max_{(x,\xi)\in \mathcal{T}} |\bar{k}(x,\xi) - \bar{k}^n(x,\xi)| & = 0.
	\end{align}
\end{subequations}
As proved in Lemma \ref{lem:ka} (relation \eqref{eq:sa}), this convergence is enabled by the 
convergence of the respective parameters $\lambda^n, \sigma^n, \theta^n, W^n, q^n$ 
(constructed in \eqref{eq:sfn}, based on the parameters associated with the solution to 
\eqref{eq:kn}, \eqref{eq:knbc}) associated with the kernels constructed in \eqref{eq:ensk} based 
on the solutions to \eqref{eq:kn}, \eqref{eq:knbc}, to the parameters 
$\lambda, \sigma, \theta, W,q$ associated with the kernels satisfying \eqref{eq:kc}, 
\eqref{eq:kcbc}. Moreover, by \eqref{eq:ucest}, the solution sequence $(u^n(t), v^n(t))$, 
constructed 
in \eqref{eq:undef} based on the solutions to \eqref{eq:n+1}, \eqref{eq:nbcuy},
 under the conditions of Theorem \ref{thm:ua}, converges (in 
$E_c$, on compact time intervals) to the solution of \eqref{eq:inf}, \eqref{eq:cbcuy}, i.e., 
\begin{equation}
	\label{eq:uninf}
	\lim_{n\to\infty}\left\|\left(\begin{smallmatrix}
		u(t) \\ v(t)
	\end{smallmatrix}\right) - \left(\begin{smallmatrix}
	u^n(t) \\ v^n(t)
	\end{smallmatrix}\right)\right\|_{E_c} = 0,
\end{equation}
for all $t \in [0,T]$, for any given $T > 0$. Note that all convergence properties stated hold in 
$L^2$ for $y$, which is 
sufficient for attaining the stated approximation result. In fact, based on the convergence 
properties of the kernels and solutions, namely, equations \eqref{eq:kainf} and \eqref{eq:uninf}, 
respectively, the control law \eqref{eq:Un} converges to \eqref{eq:Uc} 
(on compact time intervals), i.e., the following holds
\begin{align}
	\label{eq:Uinf}
	\lim_{n\to\infty} \int\limits_0^1 \left[  \frac{1}{n}\sum_{i=1}^nk^i(1,\xi)u^i(t,\xi) +
	k^{n+1}(1,\xi)v(t,\xi) \right] d\xi & = \nonumber\\
	\lim_{n\to\infty}\left\langle \left(\begin{smallmatrix}
		\mathbf{k}(1,\cdot) \\ k^{n+1}(1,\cdot)
	\end{smallmatrix}\right), \left(\begin{smallmatrix}
	\mathbf{u(t)} \\ v(t)
	\end{smallmatrix}\right) \right\rangle_E & = \nonumber \\
	\lim_{n\to\infty}\left\langle \left(\begin{smallmatrix}
		k^n(1,\cdot,\cdot) \\ \bar{k}^n(1,\cdot)
	\end{smallmatrix}\right), \left(\begin{smallmatrix}
		u^n(t) \\ v^n(t)
	\end{smallmatrix}\right) \right\rangle_{E_c} & = \nonumber \\
	\left\langle \left(\begin{smallmatrix}
		k(1,\cdot,\cdot) \\ \bar{k}(1,\cdot)
	\end{smallmatrix}\right), \left(\begin{smallmatrix}
		u(t) \\ v(t)
	\end{smallmatrix}\right) \right\rangle_{E_c} & = \nonumber \\
	\int\limits_0^1 \left[\int\limits_0^1 k(1,\xi,y)u(t,\xi,y)dy
	+ \bar{k}(1,\xi)v(t,\xi) \right]d\xi, &
\end{align}
for all $t \in [0,T]$, for any given $T > 0$.

Finally, it is interesting to note that the stability estimate 
\eqref{eq:uvdecay}, based on the previous discussion, in the limiting case $n\to\infty$ becomes, 
essentially, 
the stability estimate of the continuum system \eqref{eq:inf}, \eqref{eq:cbcuy} under the exact 
continuum control law \eqref{eq:Uc}, i.e., (see also \cite[Thm 2]{AllKrs23arxiv})
\begin{equation}
	\label{eq:uvdecayinf}
	\left\| \left( 
	\begin{smallmatrix}
		u(t) \\ v(t)
	\end{smallmatrix}
	\right) \right\|^2_{E_c} \leq \widetilde{M}^C_V\left\| \left( 
	\begin{smallmatrix}
		u(0) \\ v(0)
	\end{smallmatrix}
	\right) \right\|^2_{E_c}\exp\left(-\frac{\bar{c}_V}{M^C_V}t\right),
\end{equation}
where the coefficients $\widetilde{M}_V^C, M_V^C, \bar{c}_V$ are obtained based on the 
continuum counterparts of the respective coefficients in the 
Lyapunov-based proof of Theorem \ref{thm:stab}. This can be shown as follows. First, the 
solution $\left( (u^i(t))_{i=1}^n, v(t)\right)$ to \eqref{eq:n+1}, 
\eqref{eq:nbcuy} converges to the solution $(u(t),v(t))$ to \eqref{eq:inf}, \eqref{eq:cbcuy} 
(via the step function interpretation \eqref{eq:undef}; see Theorem~\ref{thm:ua} for details) as per 
\eqref{eq:uninf}, which results in the $E_c$ 
norms of $(u(t),v(t))$ and $(u_0,v_0)$ appearing in \eqref{eq:uvdecayinf} and corresponding to 
the $E$ norm in \eqref{eq:uvdecay}\footnote{In more detail, using \eqref{eq:undef} and 
\eqref{eq:uninf}, it holds that
\begin{equation}
\lim_{n\to\infty}\left\| \left( \begin{smallmatrix}
	\mathbf{u} \\ v
\end{smallmatrix} \right)\right\|_E = \lim_{n\to\infty}\left\| \mathcal{F}\left( \begin{smallmatrix}
\mathbf{u} \\ v
\end{smallmatrix} \right)\right\|_{E_c} = \lim_{n\to\infty}\left\| \left( \begin{smallmatrix}
u^n\\ v^n
\end{smallmatrix} \right)\right\|_{E_c} = \left\| \left( \begin{smallmatrix}
u \\ v
\end{smallmatrix} \right)\right\|_{E_c}.
\end{equation}}. Second, combining the convergence property of each of the parameter 
sequences to their respective continuum limits, which follows from \eqref{eq:sa}, together with 
the uniform essential boundedness of the parameter sequences, which follows from the step 
function-based construction, we can uniformly bound the functions appearing on the left-hand 
side of  \eqref{eq:lbd}, and hence, also bound their continuum limits. In particular, $M_\sigma$ 
and $M_W$ can be replaced by
\begin{align}
	M_\sigma^C & = \max_{x\in[0,1]} \|\sigma(x,\cdot,\cdot)\|_{L^\infty([0,1]^2;\mathbb{R})}, 
	\label{eq:Msc}\\
	M_W^C & = \max_{x\in[0,1]}\|W(x,\cdot)\|_{L^\infty([0,1];\mathbb{R})}.
\end{align}
Similarly, we can replace $M_\lambda$ by $\displaystyle M_\lambda^C = 
\max_{x\in[0,1]}\|1/\lambda(x,\cdot)\|_{L^\infty([0,1];\mathbb{R})}$.
Moreover, based on the convergence property \eqref{eq:kainf}, together with the essential 
boundedness of $(k^n, \bar{k}^n)$ (by construction), we can replace $M_k$ by
\begin{equation}
	\resizebox{.99\columnwidth}{!}{$\displaystyle 
	M_k^C = \max \left\{\max_{(x,\xi)\in \mathcal{T}}\|k(x,\xi,\cdot)\|_{L^\infty([0,1];\mathbb{R})},
	\max_{(x,\xi)\in \mathcal{T}}|\bar{k}(x,\xi)|\right\}.$}
\end{equation}
Consequently, as the inverse  kernels $(l^i)_{i=1}^{n+1}$ depend continuously on 
$(k^i)_{i=1}^{n+1}$, while the parameters $(c_{i,j})_{i,j=1}^n, (\kappa_i)_{i=1}^n$ depend 
continuously on $(k^i)_{i=1}^{n+1}$ and $(W_i)_{i=1}^n$, we can construct step function 
sequences $(l^n,\bar{l}^n), c^n, \kappa^n$ (associated with parameters $l^i, c_{i,j}, \kappa_i$, in 
a similar manner to constructions \eqref{eq:sfn}, \eqref{eq:ensk}), based on the step function 
sequences $(k^n, \bar{k}^n)$ and 
$W^n$. These sequences converge to their respective continuum 
limits $(l,\bar{l}), c, \kappa$ (in the $L^2$ sense in $y$), and they can be bounded, similarly to 
$(k^n,\bar{k}^n)$ and $W^n$. Thus, we can replace $M_\kappa$ and $M_c$ in \eqref{eq:lbd}, 
and $M_l$ in \eqref{eq:lest}, by
\begin{align}
	& M_\kappa^C = \max_{(x,\xi)\in \mathcal{T}} 
	\|\kappa(x,\xi,\cdot)\|_{L^\infty([0,1];\mathbb{R})}, \\
	& M_c^C= \max_{(x,\xi)\in \mathcal{T}}\|c(x,\xi,\cdot)\|_{L^\infty([0,1];\mathbb{R})}, \\
	& 	\resizebox{.99\columnwidth}{!}{$\displaystyle M_l^C =  \max \left\{\max_{(x,\xi)\in 
	\mathcal{T}}\|l(x,\xi,\cdot)\|_{L^\infty([0,1];\mathbb{R})}, 
	\max_{(x,\xi)\in \mathcal{T}}|\bar{l}(x,\xi)|\right\}.$} \label{eq:Mlc}
\end{align} 
Therefore, as $\delta,\epsilon\to 0$ in the limiting case 
$n \to \infty$, all involved parameters in estimate \eqref{eq:uvdecay}, namely $\widetilde{M}_V, 
c_V$, and $M_V$, can be replaced by the continuum bounds  $\widetilde{M}_V^C, 
\bar{c}_V$ and $M_V^C$, based on \eqref{eq:Msc}--\eqref{eq:Mlc}\footnote{In particular, based 
on the explicit dependence of $\widetilde{M}_V, M_V, m_V, c_V$ on 
$M_\lambda,M_W,M_k,M_\sigma,M_l,M_\kappa,M_c$, 
we can obtain the respective expressions for $\widetilde{M}_V^C, \bar{c}_V, M_V^C$, using 
\eqref{eq:Msc}--\eqref{eq:Mlc}, together with $M_\lambda^C$ and $m_\lambda^C = \displaystyle 
\max_{x\in[0,1]}\|\lambda(x,\cdot)\|_{L^\infty([0,1];\mathbb{R})}$.}, fact 
which allows us to obtain \eqref{eq:uvdecayinf}. We finally note that, by \eqref{eq:uninf} and 
\eqref{eq:Uinf}, the estimate 
\eqref{eq:uvdecayinf} holds for all $t\in[0,T]$, for any given $T > 0$, nevertheless, since 
\eqref{eq:uvdecayinf} also provides a uniform bound for the solution, the estimate 
\eqref{eq:uvdecayinf} is valid for all $t \in [0, +\infty)$.

\section{Conclusions and Discussion} \label{sec:conc}

Computation of approximate stabilizing kernels based on the continuum kernel may provide 
flexibility in computation (also in terms of the number of different kernels being computed), as 
well as it may significantly improve computational complexity (although, practically, such 
computation also depends on the sampling method chosen for the continuum kernel). This is 
confirmed in the numerical example in which computational burden of stabilizing kernels is 
significantly improved, since the approximate kernels computed based on the continuum can be 
computed in closed form, in contrast to the exact kernels that have to be computed implicitly 
based on the solution to the kernel PDEs. 

In general, we may expect that the complexity of computation of stabilizing control 
gains via the continuum approximation approach to not scale with $n$, i.e., to be 
$\mathcal{O}(1)$; 
while the complexity of computation of the exact control kernels to be 
$\mathcal{O}(n)$, i.e., to grow with the number of state components. Thus, this approach may be 
useful for computationally efficient control of large-scale PDE systems.

Another important possible usage of the approach presented may be in reducing the number of 
measurements required in a full-state feedback law, retaining closed-loop stability, in view of the 
fact that, as the number of state components increases, the sum in the control law \eqref{eq:aUn} 
essentially approximates the integral in \eqref{eq:Uc}, which in turn may be computed at a chosen 
resolution, not necessarily equal to the number of state components. Although this is also 
motivated by simulation investigations we performed using the same numerical example as in 
Section~\ref{sec:sim}, in which stability is 
preserved when we employ in \eqref{eq:aUn} a smaller number, than $n$, of terms in the sum, 
this has to be rigorously proved. We have made here the first step towards such study by formally 
establishing the connection between the solutions to 
the original $n+1$ system \eqref{eq:n+1}, \eqref{eq:nbcuy} and the solutions to the continuum 
system \eqref{eq:inf}, \eqref{eq:cbcuy}, as the number $n$ of state 
components becomes sufficiently large. This is a matter we are currently investigating.

\appendix

\setcounter{equation}{0}
\renewcommand{\theequation}{\thesubsection.\arabic{equation}}
\renewcommand{\thetheorem}{\thesubsection.\arabic{theorem}}

\subsection{Technical Propositions} \label{app:tech}

This appendix contains two technical results regarding the well-posedness of the $n+1$ system 
\eqref{eq:n+1}, \eqref{eq:nbcuy}
under Assumption~\ref{ass:n+1}, and a generic result regarding exponential stability of 
perturbed, well-posed, infinite-dimensional linear systems. Both results are direct consequences 
of existing results in the literature.
\begin{proposition}
	\label{prop:wp}
	The system \eqref{eq:n+1}, \eqref{eq:nbcuy} is well-posed on $L^2([0,1];\mathbb{R}^{n+1})$ 
	for any parameters satisfying Assumption~\ref{ass:n+1}.
	
	\begin{proof}	
		We utilize \cite[Thm 13.2.2]{JacZwaBook}, for which we need to show that the system 
		\eqref{eq:n+1}, \eqref{eq:nbcuy} can be interpreted as a port-Hamiltonian system satisfying 
		the assumptions of  this theorem. Thus, we write the system 
		\eqref{eq:n+1} in the form $z_t(t,x) = P_1(\mathcal{H}(x)z(t,x))_x + D(x) z(t,x)$, where $z = 
		(u^1, \ldots, u^n, v)^T, P_1 = \operatorname{diag}(-I_{n\times n}, 1) = P_1^T, \mathcal{H} = 
		\operatorname{diag}(\lambda_1, \ldots, \lambda_n, \mu) \in C^1([0,1]; 
		\mathbb{R}^{(n+1)\times(n+1)})$ by Assumption~\ref{ass:n+1} and $D$ is a bounded linear 
		operator. More precisely, $D = -P_1\mathcal{H}' + D_1$, where $D_1$ is such that 
		$D_1(x)z(x,t)$ corresponds to the right-hand-side of \eqref{eq:n+1}, and the the diagonal 
		operator $P_1\mathcal{H}$ satisfies the first 
		assumption of \cite[Thm 13.2.2]{JacZwaBook} by Assumption~\ref{ass:n+1}.
		
		For the boundary conditions \eqref{eq:nbcuy}, we define a vector $f = 
		(u^1(0), \ldots, u^n(0)$, $v(0), u^1(1), \ldots, u^n(1),v(1))^T$ comprising the boundary 
		values of $z$. Using this, the boundary 
		condition $u^i(\cdot,0) = q_iv(\cdot,0)$ can be written as $W_{B,1}f = 0$, where $W_{B,1} = 
		\begin{bmatrix}
			I_{n\times n} & -\mathbf{q} & 0_{n\times (n+1)}
		\end{bmatrix}$, where we denote $\mathbf{q} = (q_1, \ldots, q_n)^T$. The boundary control 
		$v(\cdot, 1) = U(\cdot)$ can be written as 
		$W_{B,2}f = U$, where $W_{B,2} = \begin{bmatrix}
			0 & \ldots & 0 & 1
		\end{bmatrix}$. Thus, the matrix $W_B = \begin{bmatrix}
			W_{B,1} \\ W_{B,2}
		\end{bmatrix}$ has full rank regardless of $\mathbf{q}$, which covers 
		the second assumption of  \cite[Thm 13.2.2]{JacZwaBook}. Thus, by  \cite[Thm 13.2.2, Rem. 
		13.2.3]{JacZwaBook}, the system \eqref{eq:n+1}, \eqref{eq:nbcuy} is well-posed.
	\end{proof}
\end{proposition}

\begin{proposition}
	\label{prop:wpstab}
	Consider a well-posed abstract Cauchy problem $\dot{z}(t) = Ax(t) + BU(t)$ on a Hilbert space 
	$Z$ under feedback control law $U(t) = (K + \Delta K)z(t)$, where $K$ and $\Delta K$ are 
	bounded linear operators. If the feedback law $U(t) = Kz(t)$ is exponentially stabilizing, then so 
	is $U(t) = Kz(t) + \Delta Kz(t)$, provided that $\|\Delta 
	K\|_{\mathcal{L}(Z,\mathbb{R})}$ is sufficiently small.
	
	\begin{proof}
		The proof relies on the well-posedness of the abstract Cauchy problem, i.e., $B$ being an 
		admissible control operator for (the strongly continuous semigroup generated by) $A$ 
		\cite[Sect. 4.2]{TucWeiBook}, and the Gearhart-Greiner-Pr\"uss Theorem \cite[Thm 
		V.1.11]{EngNagBook}. Let us denote $A_K = A+BK$ which is the generator of an exponentially 
		stable semigroup on 
		$Z$, i.e., the solution to $\dot{z}(t) = (A+BK)z(t), z(0) = z_0$ satisfies $\|z(t)\|_Z\leq 
		Me^{-\gamma t}\|z_0\|_Z$ for some $M, \gamma > 0$. Now, the system under the perturbed 
		control law can be written as $\dot{z}(t) = A_Kz(t)+B\Delta Kz(t)$, where $B$ is admissible 
		for $A_K$ by \cite[Cor. 5.5.1]{TucWeiBook}. Moreover, there exists 
		some $K_{\gamma/4} > 0$ such that \cite[Prop. 4.4.6]{TucWeiBook}
		\begin{align} 
			\label{eq:Ba}
			& \|(sI-A_K)^{-1}B\Delta K\|_{\mathcal{L}(Z)} \nonumber \\
			& \qquad  \leq \frac{K_{\gamma/4}}{\sqrt{\rm{Re}(s) + \gamma/4}}\|\Delta
			K\|_{\mathcal{L}(Z,\mathbb{R})}, \quad \forall \operatorname{Re}(s) > -\gamma/4 
			\nonumber \\
			& \qquad \leq \frac{2K_{\gamma/4}}{\sqrt{\gamma}}\|\Delta 
			K\|_{\mathcal{L}(Z,\mathbb{R})}, 
			\quad \forall 
			\operatorname{Re}(s) > 0,
		\end{align}
		where, in place of $-\gamma/4$ we could use any value on $(-\gamma, 0)$.
		
		The rest of the proof utilizes spectral arguments. Based on the identity $(sI-A_K-B\Delta K) = 
		(sI-A_K)(I + (sI -	A_K)^{-1}B\Delta K)$, we see that any point $s \in \mathbb{C}$ is 
		guaranteed to be in the resolvent set of $A_K + B\Delta K$ if it is in the resolvent set of 
		$A_K$  and $\|(sI -A_K)^{-1}B\Delta K\|_{\mathcal{L}(Z)}<1$. By the exponential stability 
		of $A_K$ and 
		\eqref{eq:Ba}, 
		these conditions can be guaranteed for all $\rm{Re}(s) > 0$ when $\|\Delta 
		K\|_{\mathcal{L}(Z,\mathbb{R})}$ is sufficiently 
		small, in addition to the resolvent operator $(sI - A_K - B\Delta K)^{-1}$ being uniformly 
		bounded for all $\rm{Re}(s) > 0$. Thus, for $\|\Delta K\|_{\mathcal{L}(Z,\mathbb{R})}$ 
		sufficiently small, we have that 
		$A_K + B\Delta K$ is exponentially stable by \cite[Thm V.1.11]{EngNagBook}, which 
		concludes the proof.
	\end{proof}
\end{proposition}
\subsection{Well-Posedness of the Continuum System} \label{app:wp}

\setcounter{equation}{0}
\setcounter{theorem}{0}

The following proposition states the well-posedness result for \eqref{eq:inf}, \eqref{eq:cbcuy}, 
analogously to what has been stated in Proposition~\ref{prop:wp} for \eqref{eq:n+1}, 
\eqref{eq:nbcuy}.
\begin{proposition}
	\label{prop:cwp}
	The continuum system \eqref{eq:inf}, \eqref{eq:cbcuy} is well-posed on $E_c = L^2([0, 
	1];L^2([0,1];\mathbb{R}))\times L^2([0,1];\mathbb{R})$ for any parameters satisfying
	Assumption~\ref{ass:inf}. That is, for any initial conditions $u_0 \in L^2([0, 
	1];L^2([0,1];\mathbb{R})), v_0 \in L^2([0,1];\mathbb{R})$ and input $U \in L_{\rm{loc}}^2([0, 
	+\infty);\mathbb{R})$, there is a unique (weak) solution to  \eqref{eq:inf}, \eqref{eq:cbcuy} 
	satisfying $(u,v) \in C([0,+\infty); E_c)$.
	
	\begin{proof}
		In order to show the well-posedness of the continuum system
		\eqref{eq:inf}, \eqref{eq:cbcuy}, due to
		\cite[Lem. 13.1.14]{JacZwaBook}, we only need to show the
		well-posedness of the pure transport part, i.e., 
		\begin{subequations}
			\label{eq:inft}%
			\begin{align}
				u_t(t,x,y) + \lambda(x,y)u_x(t,x,y)
				& = 0, \\
				v_t(t,x) - \mu(x)v_x(t,x) & = 0,
			\end{align}
		\end{subequations}
		with boundary conditions $u(t,0,y) = q(y)v(t,0), v(t,1)=U(t)$
		and initial conditions $u(0,x,y) = u_0(x,y), v(0,x) =
		v_0(x)$. For now, assuming that everything is sufficiently smooth, the
		solution of the transport equation \eqref{eq:inft} with spatially-varying transport speed
		is (cf. \cite[Prop. 3.1]{CorVaz13})
		\begin{subequations}
			\label{eq:inftsol}%
			\begin{align}
				v(t,x) & = 
				\begin{cases}
					v_0 \left(\phi_{\mu}^{-1}(\phi_{\mu}(x)+t) \right), &t -
					\int\limits_x^1 \frac{d\xi}{\mu(\xi)} \leq 0 \\
					U \left(t - \int\limits_x^1 \frac{d\xi}{\mu(\xi)}\right), & t -
					\int\limits_x^1 \frac{d\xi}{\mu(\xi)} > 0
				\end{cases}, \\
				u(t,x,y) & = 
				\begin{cases}
					u_0 \left(\phi_{\lambda_y}^{-1}(\phi_{\lambda_y}(x)-t),y\right),
					&  t - \int\limits_0^{x} \frac{d\xi}{\lambda(\xi,y)} \leq 0 \\
					q(y)v \left(t - \int\limits_0^{x} \frac{d\xi}{\lambda(\xi,y)}, 0
					\right), & t - \int\limits_0^{x} \frac{d\xi}{\lambda(\xi,y)} > 0
				\end{cases},
			\end{align}
		\end{subequations}
		where $\phi_{\mu}(x)=\int_0^x \frac{d\xi}{\mu(\xi)}$ and $\phi_{\lambda_y}(x) = \int_0^x 
		\frac{d\xi}{\lambda(\xi,y)}$, pointwise for (almost) every $y\in [0,1]$.
		The solution \eqref{eq:inftsol} is well-defined, unique, and depends continuously on $\lambda,
		\mu, v_0, u_0$ and $U$. The solution formula  (and thus, the properties of the solution) is 
		valid 
		for any $v_0 \in
		H^1([0,1]; \mathbb{R}), u_0 \in H^1([0,1];L^2([0,1];\mathbb{R})),
		$ and $U\in H^1_{\rm{loc}}([0,+\infty);\mathbb{R})$ under
		the compatibility conditions $v_0(1) = U(0)$ and
		$u_0(0,y)=q(y)v_0(0)$, i.e., \eqref{eq:inftsol} is the classical
		solution of \eqref{eq:inft}. Moreover, for any $v_0\in L^2([0,1];
		\mathbb{R}), u_0 \in L^2([0,1];L^2([0,1];\mathbb{R}))$ and $U \in
		L^2_{\rm{loc}}([0,+\infty); \mathbb{R})$, we can construct a unique weak
		solution as the limit of a sequence of classical solutions
		\cite[Sect. 2.1.3]{BasCorBook} due to $H^1$ being dense in $L^2$ (on
		one-dimensional intervals) by the Sobolev Embedding Theorem. Thus, the
		system \eqref{eq:inft} is well-posed on $L^2([0,1];
		L^2([0,1];\mathbb{R}))\times L^2([0,1];\mathbb{R})$. Finally, since
		\eqref{eq:inf}, \eqref{eq:cbcuy} only differs from \eqref{eq:inft} by
		a bounded additive perturbation (i.e., everything on the right-hand
		side of \eqref{eq:inf}), we have that \eqref{eq:inf}, \eqref{eq:cbcuy}
		is well-posed by \cite[Lem. 13.1.14]{JacZwaBook}.
	\end{proof}
\end{proposition}

\end{document}